\documentclass[aop,MSNbibl,citesort,seceqn,dvips]{arximspdf}

\doi{10.
/09-AOP459}
\volume{37}
\issue{6}
\pubyear{2009}
\firstpage{2093}
\lastpage{2134}

\makeatletter

\def\iint{\int\!\!\int}

\newtheorem{theorem}{Theorem}[section]
\newproclaim{remark}{Remark}
\newtheorem{proposition}[remark]{Proposition}
\newtheorem{lemma}[remark]{Lemma}

\newproclaim{step}{Step}

\makeatother

\begin{document}
\begin{frontmatter}

\title{Variations and estimators for self-similarity parameters via
Malliavin calculus}
\runtitle{Variations and estimators via Malliavin calculus}

\begin{aug}
\author[A]{\fnms{Ciprian A.} \snm{Tudor}\corref{}\ead[label=e1]{tudor@univ-paris1.fr}} and
\author[B]{\fnms{Frederi G.} \snm{Viens}\thanksref{t1}\ead[label=e2]{viens@stat.purdue.edu}}
\runauthor{C. A. Tudor and F. G. Viens}
\affiliation{University of Paris 1 and Purdue University}
\address[A]{SAMOS-MATISSE\\
Centre d'Economie de La Sorbonne\\
Universit\'{e} de Paris 1 Panth\'eon-Sorbonne\\
90, rue de Tolbiac\\
75634, Paris\\
France\\
\printead{e1}}
\address[B]{Department of Statistics\\
\quad and Department of Mathematics\\
Purdue University \\
150 N. University Street\\
West Lafayette, Illinois 47907-2067\\
USA\\
\printead{e2}}
\end{aug}
\thankstext{t1}{Supported in part by NSF Grant 06-06-615.}

\received{\smonth{6} \syear{2008}}
\revised{\smonth{11} \syear{2008}}

%
\begin{abstract}
Using multiple stochastic integrals and the Malliavin calculus, we
analyze the
asymptotic behavior of quadratic variations for a specific non-Gaussian
self-similar process, the Rosenblatt process. We apply our results to the
design of strongly consistent statistical estimators for the self-similarity
parameter $H$. Although, in the case of the Rosenblatt process, our estimator
has non-Gaussian asymptotics for all $H>1/2$, we show the remarkable fact
that the process's data at time $1$ can be used to construct a distinct,
compensated estimator with Gaussian asymptotics for $H\in(1/2,2/3)$.
\end{abstract}

%
\begin{keyword}[class=AMS]
\kwd[Primary ]{60F05}
\kwd{60H05}
\kwd[; secondary ]{60G18}
\kwd{62F12}.
\end{keyword}
\begin{keyword}
\kwd{Multiple stochastic integral}
\kwd{Hermite process}
\kwd{fractional Brownian motion}
\kwd{Rosenblatt process}
\kwd{Malliavin calculus}
\kwd{noncentral limit theorem}
\kwd{quadratic variation}
\kwd{Hurst parameter}
\kwd{self-similarity}
\kwd{statistical estimation}.
\end{keyword}

\end{frontmatter}

\section{Introduction}

\subsection{Context and motivation}

A \textit{self-similar process} is a stochastic process such that any
part of
its trajectory is invariant under time scaling. Self-similar
processes are of considerable interest in practice in modeling
various phenomena, including internet traffic (see, e.g.,
\cite{WiTaLeWi}), hydrology (see, e.g., \cite{Leod}) or economics
(see, e.g., \cite{Mand,WiTaTe}). In various applications,
empirical data also shows strong correlation of observations,
indicating the presence, in addition to self-similarity, of
long-range dependence. We refer to the monographs \cite{EM} and
\cite{SaTa} for various properties of, and fields of application for,
such processes.

The motivation for this work is to examine non-Gaussian self-similar
processes using tools from stochastic analysis. We will focus our
attention on a special process of this type, the so-called \textit{Rosenblatt
process}. This belongs to a class of self-similar processes which also
exhibit long-range dependence and which appear as limits in the
so-called noncentral limit theorem: the class of \textit{Hermite}
processes. We study the behavior of the quadratic variations for the
Rosenblatt process $Z$, which is related to recent results by
\cite{No,NoNu,LeLu}, and we apply the results to the
study of estimators for the self-similarity parameter of $Z$.
Recently, results on variations or weighted quadratic variations of
fractional Brownian motion were obtained in \cite{No,NoNu,LeLu},
among others. The Hermite processes were
introduced by Taqqu (see \cite{Ta1} and \cite{Ta2}) and by Dobrushin
and Major (see \cite{DM}). The Hermite process of order $q\geq1$ can
be written, for every $t\geq0$, as%
%
\begin{eqnarray}\label{hermite}\hspace*{30pt}
Z_{H}^{q}(t)&=&c(H,q)\nonumber\\[-8pt]\\[-8pt]
&&{}\times \int_{\mathbb{R}^{q}} \Biggl[ \int_{0}^{t} \Biggl(
\prod_{i=1}^{q}(s-y_{i})_{+}^{- ( {1}/{2}+({1-H})/{q} )
} \Biggr) \,ds \Biggr]\,dW(y_{1})\cdots dW(y_{q}),\nonumber
\end{eqnarray}
where $c(H,q)$ is an explicit positive constant depending on $q$ and
$H$, and
such that $\mathbf{E} ( Z_{H}^{q}(1)^{2} ) =1$, $x_{+}=\max(x,0)$,
the self-similarity (Hurst) parameter $H$ belongs to the interval
$(\frac{1}%
{2},1)$ and the above integral is a multiple Wiener--It\^{o}
stochastic integral with respect to a two-sided Brownian motion
$ ( W(y) ) _{y\in\mathbb{R}}$ (see \cite{Nbook}). We
note that the Hermite processes of order $q>1$, which are
non-Gaussian, have only been defined for $H>\frac{1}{2}$; how to
define these processes for $H\leq\frac{1}{2}$ is still an open
problem.

The case $q=1$ is the well-known fractional Brownian motion (fBm): this is
Gaussian. One recognizes that when $q=1$, (\ref{hermite}) is the moving
average representation of fractional Brownian motion. The Rosenblatt process
is the case $q=2$. All Hermite processes share the following basic properties:

\begin{itemize}
\item they exhibit long-range dependence (the long-range covariance
decays at
the rate of the nonsummable power function $n^{2H-2}$);
\item they are $H$-self-similar, in the sense that for any $c>0$, $ (
Z_{H}^{q}(ct) ) _{t\geq0}$ and $ ( c^{H} 
Z_{H}^{q}(t) )
_{t\geq0}$ are equal in distribution;
\item they have stationary increments, that is, the distribution of $ (
Z_{H}^{q}(t+h)-Z_{H}^{q}(h) ) _{t\geq0}$ does not depend on $h>0$;
\item they share the same covariance function,%
\[
\mathbf{E} [ Z_{H}^{q}(t)Z_{H}^{q}(s) ] =:R^{H}(t,s)=\tfrac{1}%
{2} ( t^{2H}+s^{2H}-|t-s|^{2H} ),\qquad s,t\geq0,
\]
so, for every $s,t\geq0$, the expected squared increment of the
Hermite process is%
%
\begin{equation}\label{zt-zs}%
\mathbf{E} \bigl[ \bigl( Z_{H}^{q}(t)-Z_{H}^{q}(s) \bigr) ^{2} \bigr]
=|t-s|^{2H}
\end{equation}
from which it follows by Kolmogorov's continuity criterion, and the
fact that
each $L^{p} ( \Omega) $-norm of the increment of $Z_{H}^{q}$ over
$[s,t]$ is commensurate with its $L^{2} ( \Omega) $-norm, that
this process is almost surely H\"{o}ld{e}r continuous of any order
$\delta<H$;
\item the $q$th Hermite process lives in the so-called $q$th Wiener
chaos of
the underlying Wiener process $W$ since it is a $q$th order Wiener integral.
\end{itemize}

The stochastic analysis of fBm has been developed intensively in recent years
and its applications are numerous. Other Hermite processes are less
well studied, but
are still of interest due to their long-range dependence, self-similarity
and stationarity of increments. The great popularity of fBm in modeling
is due
to these properties and fBm is preferred over higher order Hermite
processes because it is a Gaussian process and because its calculus is
much easier. However,
in concrete situations, when empirical data attests to the presence of
self-similarity and long memory without the Gaussian property, one can
use a
Hermite process living in a higher chaos.

The Hurst parameter $H$ characterizes all of the important properties
of a
Hermite process, as seen above. Therefore, properly estimating $H$ is
of the
utmost importance. Several statistics have been introduced to this end, such
as wavelets, $k$-variations, variograms, maximum likelihood estimators and
spectral methods. Information on these various approaches can be found
in the
book by Beran \cite{B}.

In this paper, we will use variation statistics to estimate $H$. Let us recall
the context. Suppose that a process $(X_{t})_{t\in[0,1]}$ is
observed at
discrete times $\{0,\frac{1}{N},\ldots,\frac{N-1}{N},1\}$ and let $a$
be a
``filter'' of length $l\geq0$ and $p\geq1$ a
fixed power; that is, $a$ is an $l+1$-dimensional vector $a=(a_{0}%
,a_{1},\ldots,a_{l})$ such that $\sum_{q=0}^{l}a_{q}q^{r}=0$ for $0\leq
r\leq
p-1$ and $\sum_{q=0}^{l}a_{q}q^{p}\not=0$. The $k$-variation statistic
associated to the filter $a$ is then defined as
\[
V_{N}(k,a)=\frac{1}{N-l}\sum_{i=l}^{N-1} \biggl[ \frac{ \vert
V_{a} ( {i}/{N} ) \vert^{k}}{\mathbf{E} [
| V_{a} ( {i}/{N} ) | ^{k} ]
}-1 \biggr],
\]
where, for $i\in\{l,\ldots,N\}$,%
\[
V_{a} \biggl( \frac{i}{N} \biggr) =\sum_{q=0}^{l}a_{q}X \biggl( \frac{i-q}%
{N} \biggr).
\]
When $X$ is fBm, these statistics are used to derive strongly consistent
estimators for the Hurst parameter and their associated normal convergence
results. A detailed study can be found in \cite{GuLe} and \cite{LaIs}
or, more
recently, in \cite{coeur}. The behavior of $V_{N}(k,a)$ is used to derive
similar behaviors for the corresponding estimators. The basic result
for fBm
is that, if $p>H+\frac{1}{4}$, then the renormalized $k$-variation
$V_{N}(k,a)$ converges to a standard normal distribution. The easiest
and most
natural case is that of the filter $a=\{1,-1\}$, in which case $p=1$;
one then
has the restriction $H<\frac{3}{4}$. The techniques used to prove such
convergence in the fBm case in the above references are strongly
related to
the Gaussian property of the observations; they appear not to extend to
non-Gaussian situations.

Our purpose here is to develop new techniques that can be applied to
both the
fBm case and to other non-Gaussian self-similar processes. Since this
is the first
attempt in such a direction, we keep things as simple as possible: we treat
the case of the filter $a=\{1,-1\}$ with a $k$-variation order $=2$ (quadratic
variation), but the method can be generalized. As announced above, we further
specialize to the simplest non-Gaussian Hermite process, that is, the
one of order
$2$, the Rosenblatt process. We now give a short overview of our
results (a
more detailed summary of these facts is given in the next subsection). We
obtain that, after suitable normalization, the quadratic variation statistic
of the Rosenblatt process converges to a Rosenblatt random variable
with the
same self-similarity order; in fact, this random variable is the
observed value
of the original Rosenblatt process at time $1$ and the convergence
occurs in
the mean square. More precisely, the quadratic variation statistic can be
decomposed into the sum of two terms: a term in the fourth Wiener chaos
(i.e., an iterated integral of order 4 with respect to the Wiener
process) and a
term in the second Wiener chaos. The fourth Wiener chaos term is well behaved,
in the sense that it has a Gaussian limit in distribution, but the second
Wiener chaos term is ill behaved, in the sense that its asymptotics are
non-Gaussian and are, in fact, Rosenblatt-distributed. This term, being
of a
higher order than the well-behaved one, is responsible for the asymptotics
of the entire statistic. But, since its convergence occurs in the mean-square
and the limit is observed, we can construct an adjusted variation by
subtracting the contribution of the ill-behaved term. We find an
estimator for
the self-similarity parameter of the Rosenblatt process, based on observed
data, whose asymptotic distribution is normal.

Our main tools are the Malliavin calculus, the Wiener--It\^{o} chaos
expansions and recent results on the convergence of multiple
stochastic integrals proved in \cite{HN,NOT,NP} and
\cite{PT}. The key point is the following: if the observed process
$X$ lives in some Wiener chaos of finite order, then the statistic
$V_{N}$ can be decomposed, using product formulas and Wiener chaos
calculus, into a finite sum of multiple integrals. One can then
attempt to apply the criteria in \cite{NOT} to study the convergence
in law of such sequences and to derive asymptotic normality results (or
to demonstrate the lack thereof) on the estimators for the Hurst
parameter of
the observed process. The criteria in \cite{NOT} are necessary and
sufficient conditions for convergence to the Gaussian law; in some
instances, these criteria fail (e.g., the fBm case with $H>3/4$), in
which case, a proof of nonnormal convergence ``by
hand,'' working directly with the chaoses, can be
employed. It is the basic Wiener chaos calculus that makes this
possible.

\subsection{Summary of results}

We now summarize the main results of this paper in some detail. As stated
above, we use quadratic variation with $a=\{1,-1\}$. We consider the
two following
processes, observed at the discrete times $ \{ i/N \}
_{i=0}^{N}$: the fBm process $X=B$ and the Rosenblatt process $X=Z$. In
either case, the standardized quadratic variation and the Hurst parameter
estimator are given, respectively, by%
%
\begin{eqnarray}\qquad
\label{VNdef0}
V_{N} &=& V_{N} ( 2, \{ -1,1 \} ) :=\frac{1}{N}%
\sum_{i=1}^{N} \biggl( \frac{ \vert X ( i/N ) -X ( (
i-1 ) /N ) \vert^{2}}{N^{-2H}}-1 \biggr) ,\\
\label{HN}%
\hat{H}_{N} &=& \hat{H}_{N}(2,\{-1,1\}):=\frac{1}{2}-\frac{1}{2\log
N}\log
\sum_{i=1}^{N} \biggl( X\biggl(\frac{i}{N}\biggr)-X\biggl(\frac{i-1}{N}\biggr) \biggr) ^{2}.
\end{eqnarray}
We choose to use the normalization $\frac{1}{N}$ in the definition of $V_{N}$
(as, e.g., in \cite{coeur}), although, in the literature, it sometimes
does not
appear. The  $H$-dependent constants $c_{j,H}$ (et al.) referred to
below are
defined explicitly in (\ref{dh}), (\ref{ah}), (\ref{c3H}), (\ref{c1H}), (\ref{e1H})
and (\ref{f1H}). Here, and throughout, $L^{2} (
\Omega) $ denotes the set of square-integrable random variables
measurable with respect to the sigma field generated by $W$. This
sigma-field is the
same as that generated by $B$ or by $Z$. The term ``Rosenblatt
random variable'' denotes a random variable whose distribution is the
same as that of $Z ( 1 ) $.

We first recall the followings facts, relative to fractional
Brownian motion:
\begin{enumerate}
\item if $X=B$ and $H\in(1/2,3/4)$, then:
\begin{enumerate}[(a)]
\item[(a)] $\sqrt{N/c_{1,H}}V_{N}$ converges in distribution to the standard
normal law;
\item[(b)] $\sqrt{N}\log(N)\frac{2}{\sqrt{c_{1,H}}} ( \hat{H}_{N}-H ) $
converges in distribution to the standard normal law;
\end{enumerate}
\item if $X=B$ and $H\in(3/4,1)$, then:
\begin{enumerate}[(a)]
\item[(a)] $\sqrt{N^{4-4H}/c_{2,H}}V_{N}$ converges in $L^{2} ( \Omega)
$ to a standard Rosenblatt random variable with parameter $H_{0}=2H-1$;
\item[(b)] $N^{1-H}\log(N)\frac{2}{\sqrt{c_{2,H}}} ( \hat{H}_{N}-H ) $
converges in $L^{2} ( \Omega) $ to the same standard Rosenblatt
random variable;
\end{enumerate}
\item if $X=B$ and $H=3/4$, then:
\begin{enumerate}[(a)]
\item[(a)] $\sqrt{N/ ( c_{1,H}^{\prime}\log N ) }V_{N}$ converges in
distribution to the standard normal law;
\item[(b)] $\sqrt{N\log N}\frac{2}{\sqrt{c_{1,H}^{\prime}}}(\hat{H}_{N}(2,a)-H)$
converges in distribution to the standard normal law.
\end{enumerate}

The convergences for the standardized $V_{N}$'s in points 1(a) and
2(a) have been known for some time, in works such as \cite{Ta2} and
\cite{Ha}. Lately, even stronger results, which also give error
bounds, have been proven. We refer to \cite{NoPe2} for the
one-dimensional case and $H\in(0,\frac{3}{4})$, \cite{BN} for then
one-dimensional case and $H\in[\frac{3}{4},1)$ and to
\cite{NPR} for the multidimensional case and $H\in(0,\frac{3}{4})$.

In this paper, we prove the following results for the Rosenblatt
process $X=Z$
as $N\rightarrow\infty$:
\item if $X=Z$ and $H\in( 1/2,1 ) $, then with $c_{3,H}$ in
(\ref{c3H}),
\begin{enumerate}[(a)]
\item[(a)] $N^{1-H}V_{N}(2,a)/ ( c_{3,H} ) $ converges in $L^{2} (
\Omega) $ to the Rosenblatt random variable $Z ( 1 ) $;
\item[(b)] $\frac{N^{1-H}}{2c_{3,H}}\log( N )  (\hat{H}_{N}(2,a)-H)$
converges in $L^{2} ( \Omega) $ to the same Rosenblatt random
variable $Z ( 1 ) $;
\end{enumerate}
\item if $X=Z$ and $H\in( 1/2,2/3 ) $, then, with $e_{1,H}$ and
$f_{1,H}$ in (\ref{e1H}) and (\ref{f1H}),
\begin{enumerate}[(a)]
\item[(a)] $\frac{\sqrt{N}}{\sqrt{e_{1,H}+f_{1,H}}} [ V_{N}(2,a)-\frac
{\sqrt{c_{3,H}}}{N^{1-H}}Z(1) ] $ converges in distribution to the
standard normal law;
\item[(b)] $\frac{\sqrt{N}}{\sqrt{e_{1,H}+f_{1,H}}} [ 2\log( N )
(H-\hat{H}_{N}(2,a))-\frac{\sqrt{c_{3,H}}}{N^{1-H}}Z(1) ] $ converges
in distribution to the standard normal law.
\end{enumerate}
\end{enumerate}

Note that $Z ( 1 ) $ is the \textit{actual observed value} of the
Rosenblatt process at time $1$, which is why it is legitimate to
include it in
a formula for an estimator. Points~4 and~5 are new results. The subject of
variations and statistics for the Rosenblatt process has thus far
received too narrow a
treatment in the literature, presumably because standard techniques inherited
from the noncentral limit theorem (and sometimes based on the Fourier
transform formula for the driving Gaussian process) are difficult to apply
(see \cite{BrMa,DM,Ta2}). Our Wiener chaos calculus approach
allows us to show that the standardized quadratic variation and corresponding
estimator both converge to a Rosenblatt random variable in $L^{2} (
\Omega) $. Here, our method has a crucial advantage: we are able to
determine which Rosenblatt random variable it converges to: it is none other
than the observed value $Z ( 1 ) $. The fact that we are able to prove
$L^{2} ( \Omega) $ convergence, not just convergence in
distribution, is crucial. Indeed, when $H<2/3$, subtracting an appropriately
normalized version of this observed value from the quadratic variation
and its
associated estimator, we prove that asymptotic normality does hold in this
case. This unexpected result has important consequences for the
statistics of
the Rosenblatt process since it permits the use of standard techniques in
parameter estimation and testing.\looseness=-1

Our asymptotic normality result for the Rosenblatt process was
specifically made possible by showing that $V_{N}$ can be
decomposed into two terms: a term $T_{4}$ in the fourth Wiener chaos
and a term $T_{2}$ in the second Wiener chaos. While the
second-Wiener-chaos term $T_{2}$ always converges to the Rosenblatt
random variable $Z ( 1 ) $, the fourth chaos term $T_{4}$ converges
to a Gaussian random variable for $H\leq3/4$. We conjecture that this
asymptotic normality should also occur for Hermite processes of
higher order $q\geq3$ and that the threshold $H=3/4$ is universal.
The threshold $H<2/3$ in the results above comes from the
discrepancy that exists between a normalized $T_{2}$ and its
observed limit $Z ( 1 ) $. If we were to rephrase results~4
and 5 above, with $T_{2}$ instead of $Z ( 1 ) $ (which
is not a legitimate operation when defining an estimator since
$T_{2}$ is not observed), the threshold would be $H\leq3/4$ and the
constant $f_{1,H}$ would vanish.

Beyond our basic interest concerning parameter estimation problems, let us
situate our paper in the context of some recent and interesting works
on the
asymptotic behavior of $p$-variations (or weighted variations) for Gaussian
processes, namely the papers \cite{LeLu,MaRo,No,NoNu} and
\cite{Swa}. These recent papers study the behavior of sequences of the
type%
\[
\sum_{i=1}^{N}h \bigl( X \bigl( (i-1)/N \bigr) \bigr) \biggl(
\frac{ \vert X ( i/N ) -X ( ( i-1 ) /N )
\vert^{2}}{N^{-2H}}-1 \biggr),
\]
where $X$ is a Gaussian process (fractional Brownian motion in
\cite{LeLu,No} and \cite{NoNu}, and the solution of the heat
equation driven by a space-time white noise in \cite{Swa}) or the
iterated Brownian motion in \cite{NoPe} and $h$ is a regular
deterministic function. In the fractional Brownian motion case, the
behavior of such sums varies according to the values of the Hurst
parameter, the limit sometimes being a conditionally Gaussian random
variable, sometimes a deterministic Riemann\vadjust{\goodbreak} integral and sometimes a
pathwise integral with respect to a Hermite process. We believe that our
work is the first to tackle a non-Gaussian case, that is, when the
process $X$ above is a Rosenblatt process. Although we restrict
ourselves to the case when $h\equiv1$, we still observe the
appearance of interesting limits, depending on the Hurst parameter:
while, in general, the limit of the suitably normalized sequence is a
Rosenblatt random variable (with the same Hurst parameter $H$ as the
data, which poses a slight problem for statistical applications),
the adjusted variations (i.e., the sequences obtained by subtracting
precisely the portion responsible for the non-Gaussian
convergence) do converge to a Gaussian limit for $H\in(1/2,2/3)$.

This article is structured as follows. Section \ref{Prelim} presents
preliminaries on fractional stochastic analysis. Section \ref{nonGauss}
contains proofs of our results for the non-Gaussian Rosenblatt process. Some
calculations are recorded as lemmas that are proven in the \hyperref
[Append]{Appendix}.
Section \ref{Stat} establishes our parameter estimation
results, which follow almost trivially from the theorems in Section~\ref{nonGauss}.

\section{Preliminaries}\label{Prelim}

Here, we describe the elements from stochastic analysis that we will
need in
the paper. Consider ${\mathcal{H}}$, a real, separable Hilbert space
and $(B
(\varphi), \varphi\in{\mathcal{H}})$, an isonormal Gaussian process,
that is, a
centered Gaussian family of random variables such that $\mathbf{E} (
B(\varphi) B(\psi) ) = \langle\varphi, \psi\rangle_{{\mathcal{H}}}$.

Denote by $I_{n}$ the multiple stochastic integral with respect to $B$ (see
\cite{Nbook} and \cite{Ustu}). This $I_{n}$ is actually an isometry between the Hilbert space
${\mathcal{H}}^{\odot n}$ (symmetric tensor product) equipped with the scaled
norm $\frac{1}{\sqrt{n!}}\Vert\cdot\Vert_{{\mathcal{H}}^{\otimes n}}$
and the
Wiener chaos of order~$n$ which is defined as the closed linear span of the
random variables $H_{n}(B(\varphi))$, where $\varphi\in{\mathcal{H}},
\Vert\varphi\Vert_{{\mathcal{H}}}=1$ and $H_{n}$ is the Hermite
polynomial of
degree~$n$.

We recall that any square-integrable random variable which is
measurable with
respect to the $\sigma$-algebra generated by $B$ can be expanded into an
orthogonal sum of multiple stochastic integrals,
\[
F=\sum_{n\geq0}I_{n}(f_{n}),
\]
where $f_{n}\in{\mathcal{H}}^{\odot n}$ are (uniquely determined) symmetric
functions and $I_{0}(f_{0})=\mathbf{E} [ F ] $.

In this paper, we actually use only multiple integrals with respect
to the standard Wiener process with time horizon $[0,1]$ and, in this
case, we will always have ${\mathcal{H}}= L^{2} ( [0,1] ) $.
This notation will be used throughout the paper.

We will need the general formula for calculating products of Wiener
chaos integrals of any orders, $p$ and $q$, for any symmetric
integrands $f\in
\mathcal{H}^{\odot p}$ and $g\in\mathcal{H}^{\odot q}$; it is%
%
\begin{equation}\label{product}%
I_{p}(f)I_{q}(g)=\sum_{r=0}^{p\wedge q}r!\pmatrix{p\cr r}\pmatrix{q\cr
r}%
I_{p+q-2r}(f\otimes_{r}g),
\end{equation}
as given, for instance, in Nualart's book \cite{Nbook}, Proposition 1.1.3;
the contraction $f\otimes_{r}g$ is the element of ${\mathcal{H}}%
^{\otimes(p+q-2r)}$ defined by
%
\begin{eqnarray} \label{contra}%
&& (f\otimes_{r} g) ( s_{1}, \ldots, s_{p-r}, t_{1}, \ldots,
t_{q-r})\nonumber\\
&&\qquad =\int_{[0,T] ^{p+q-2r} } f( s_{1}, \ldots, s_{p-r}, u_{1},
\ldots,u_{r})\\
&&\qquad\quad\hspace*{48.1pt}{}\times g(t_{1}, \ldots, t_{q-r},u_{1}, \ldots, u_{r})
\,du_{1}\cdots du_{r} .\nonumber
\end{eqnarray}
We now introduce the Malliavin derivative for random variables in a
chaos of finite order. If $f\in{\mathcal{H}}^{\odot n}$, we will use
the following rule to differentiate in the Malliavin sense:
\[
D_{t}I_{n}(f)=nI_{n-1}(f_{n}(\cdot,t)),\qquad t\in[0,1].
\]

It is possible to characterize the convergence in distribution of a sequence
of multiple integrals to the standard normal law. We will use the following
result (see Theorem 4 in \cite{NOT}, also \cite{NP}).
\begin{theorem}
\label{NOTcrit}Fix $n\geq2$ and let $(F_{k},k\geq1)$, $F_{k}=I_{n}(f_{k})$
(with $f_{k}\in{\mathcal{H}}^{\odot n}$ for every $k\geq1$), be a
sequence of
square-integrable random variables in the $n$th Wiener chaos such that
$\mathbf{E}[F_{k}^{2}]\rightarrow1$ as $k\rightarrow\infty.$ The
following are then equivalent:
\begin{longlist}
\item the sequence $(F_{k}) _{k\geq0} $ converges in distribution to the
normal law $\mathcal{{N}} (0,1)$;
\item$\mathbf{E}[F_{k}^{4}]\rightarrow3$ as $k\rightarrow
\infty$;
\eject
\item for all $1\leq l\leq n-1$, it holds that $\lim_{k\rightarrow
\infty}\Vert f_{k}\otimes_{l}f_{k}\Vert_{{\mathcal{H}}^{\otimes2(n-l)}}=0$;
\item$\Vert DF_{k}\Vert_{{\mathcal{H}}}^{2}\rightarrow n$ in
$L^{2}(\Omega)$ as $k\rightarrow\infty$, where $D$ is the Malliavin derivative
with respect to $B$.
\end{longlist}
\end{theorem}

Criterion (iv) is due to \cite{NOT}; we will refer to it as the
\textit{Nualart--Ortiz-Latorre criterion}. A multidimensional version
of the above
theorem has been proven in \cite{PT} (see also \cite{NOT}).

\section{Variations for the Rosenblatt process}\label{nonGauss}

Our\vspace*{1pt} observed process is a Rosenblatt process $(Z(t))_{t\in[0,1]}$ with
self-similarity parameter $H\in(\frac{1}{2},1)$. This centered process is
self-similar with stationary increments and lives in the second Wiener chaos.
Its covariance is identical to that of fractional Brownian motion. Our
goal is to estimate its self-similarity parameter $H$ from discrete
observations of its sample paths. As far as we know, this direction has seen
little or no attention in the literature and the classical techniques (e.g.,
the ones from \cite{DM,Ta1} and \cite{Ta2}) do not work well for it.
Therefore, the use of the Malliavin calculus and multiple stochastic integrals
is of interest.\vadjust{\goodbreak}

The Rosenblatt process can be represented as follows (see \cite{Tud}): for
every $t\in[0,1]$,
%
\begin{eqnarray}\label{rose2}\quad
Z^{H}(t)&:=& Z(t)\nonumber\\
&=& d(H)\int_{0}^{t}\int_{0}^{t} \biggl[ \int_{y_{1}\vee y_{2}}%
^{t}\partial_{1}K^{H^{\prime}}(u,y_{1})\\
&&\hspace*{51.4pt}\hspace*{34.3pt}{}\times \partial_{1}K^{H^{\prime}}%
(u,y_{2})\,du \biggr]\,dW(y_{1})\,dW(y_{2}),\nonumber
\end{eqnarray}
where $(W(t),t\in[0,1])$ is some standard Brownian motion,
$K^{H^{\prime}}$ is the standard kernel of fractional Brownian
motion of index $H^{\prime}$ (see any
reference on fBm, such as \cite{Nbook}, Chapter 5) and%
%
\begin{equation}\label{dh}%
H^{\prime}=\frac{H+1}{2} \quad\mbox{and}\quad d(H)=\frac{ ( 2(2H-1) )
^{1/2}}{ ( H+1 ) H^{1/2}}.
\end{equation}

For every $t\in[0,1]$, we will denote the kernel of the Rosenblatt
process with respect to $W$ by
%
\begin{eqnarray}\label{lt}\qquad\quad
L_{t}^{H}(y_{1},y_{2})
&:=&
L_{t}(y_{1},y_{2})
\nonumber\\[-8pt]\\[-8pt]
&:=&
d(H) \biggl[ \int_{y_{1}\vee y_{2}
}^{t}\partial_{1}K^{H^{\prime}}(u,y_{1})\,\partial_{1}K^{H^{\prime}}%
(u,y_{2})\,du \biggr] 1_{[0,t]^{2}}(y_{1},y_{2}).\nonumber
\end{eqnarray}
In other words, in particular, for every $t$,
\[
Z(t)=I_{2} ( L_{t}(\cdot) ),
\]
where $I_{2}$ denotes the multiple integral of order 2 introduced in Section
\ref{Prelim}.

Now, consider the filter $a=\{-1,1\}$ and the $2$-variations given by
\begin{eqnarray*}
V_{N}(2,a) &=& \frac{1}{N}\sum_{i=1}^{N}\frac{ ( Z({i}/{N}%
)-Z(({i-1})/{N}) ) ^{2}}{\mathbf{E} ( Z({i}/{N})-Z(({i-1})/{N}) ) ^{2}}-1 \\
&=& N^{2H-1}\sum_{i=1}^{N} \biggl[ \biggl( Z\biggl(\frac{i}%
{N}\biggr)-Z\biggl(\frac{i-1}{N}\biggr) \biggr) ^{2}-N^{-2H} \biggr] .
\end{eqnarray*}
The product formula for multiple Wiener--It\^{o} integrals (\ref{product})
yields%
\[
I_{2}(f)^{2}=I_{4}(f\otimes f)+4I_{2}(f\otimes_{1}f)+2\Vert f\Vert
^{2}_{ L^{2}([0,1]^{2})}.
\]
Setting, for $i=1, \ldots,N$,
%
\begin{equation}\label{AiRos}%
A_{i}:=L_{{i}/{N}}-L_{({i-1})/{N}},
\end{equation}
we can thus write
\[
\biggl( Z\biggl(\frac{i}{N}\biggr)-Z\biggl(\frac{i-1}{N}\biggr) \biggr) ^{2}= ( I_{2}%
(A_{i}) ) ^{2}=I_{4}(A_{i}\otimes A_{i})+4I_{2}(A_{i}\otimes_{1}%
A_{i})+N^{-2H}
\]
and this implies that the 2-variation is decomposed into a fourth chaos term
and a second chaos term:
\[
V_{N}(2,a)=N^{2H-1}\sum_{i=1}^{N} \bigl( I_{4}(A_{i}\otimes A_{i}%
)+4I_{2}(A_{i}\otimes_{1}A_{i}) \bigr) :=T_{4}+T_{2}.
\]
A detailed study of the two terms above will shed light on some interesting
facts: if $H\leq\frac{3}{4}$, then the term $T_{4}$ continues to exhibit
``normal'' behavior (when renormalized, it
converges in law to a Gaussian distribution), while the term $T_{2}$, which
turns out to be dominant, never converges to a Gaussian law. One can
say that
the second Wiener chaos portion is ``ill behaved''; however, once it is
subtracted, one obtains a
sequence converging to $\mathcal{N}(0,1)$ for $H\in(\frac{1}{2}, \frac
{2}%
{3})$, which has an impact on statistical applications.

\subsection{Expectation evaluations}

\subsubsection{The term $T_{2}$}

Let us evaluate the mean square of the second term,
\[
T_{2}:=4N^{2H-1}\sum_{i=1}^{N}I_{2}(A_{i}\otimes_{1}A_{i}).
\]
We use the notation $I_{i}= ( \frac{i-1}{N},\frac{i}{N} ] $ for
$i=1,\ldots,N$. The contraction $A_{i}\otimes_{1}A_{i}$ is given by
%
\begin{eqnarray}\label{aieaieaie}%
&&(A_{i}\otimes_{1}A_{i})(y_{1},y_{2})
\nonumber\\
&&\qquad= \int_{0}^{1}A_{i}(x,y_{1})A_{i}(x,y_{2})\,dx\nonumber\\
&&\qquad= d(H)^{2}\int_{0}^{1}dx\, 1_{[0,{i}/{N}]}(y_{1}\vee
x)1_{[0,{i}/{N}]}(y_{2}\vee x)\nonumber\\
&&\qquad\quad\hspace*{27.7pt}{}\times \biggl( \int_{x\vee y_{1}}^{{i}/{N}}\partial_{1}K^{H^{\prime}%
}(u,x)\,\partial_{1}K^{H^{\prime}}(u,y_{1})\,du\nonumber\\
&&\hspace*{78.16pt}{} - 1_{[0,({i-1})/{N}]}(y_{1}\vee
x)\\
&&\hspace*{88.4pt}{}\times\int_{x\vee y_{1}}^{({i-1})/{N}}\partial_{1}K^{H^{\prime
}}(u,x)\,\partial
_{1}K^{H^{\prime}}(u,y_{1})\,du \biggr) \nonumber\\
&&\qquad\quad\hspace*{27.7pt}{}\times  \biggl( \int_{x\vee y_{2}}^{{i}/{N}}\partial_{1}K^{H^{\prime}%
}(v,x)\,\partial_{1}K^{H^{\prime}}(v,y_{2})\,dv\nonumber\\
&&\hspace*{78.16pt}{}-1_{[0,
({i-1})/{N}]}(y_{2}\vee
x)\nonumber\\
&&\hspace*{88.4pt}{}\times\int_{x\vee y_{2}}^{({i-1})/{N}}\partial_{1}K^{H^{\prime
}}(v,x)\,\partial
_{1}K^{H^{\prime}}(v,y_{2})\,dv \biggr) .\nonumber
\end{eqnarray}
Defining
%
\begin{equation}\label{ah}%
a ( H ) :=H^{\prime} ( 2H^{\prime}-1 ) =H (
H+1 ) /2,
\end{equation}
note the following fact (see \cite{Nbook}, Chapter 5):%
%
\begin{equation}\label{souvent}%
\int_{0}^{u\wedge v}\partial_{1}K^{H^{\prime}}(u,y_{1})\,\partial_{1}%
K^{H^{\prime}}(v,y_{1})\,dy_{1}=a(H)|u-v|^{2H^{\prime}-2};
\end{equation}
in fact, this relation can easily be derived from $\int_{0}^{u\wedge
v}K^{H^{\prime}}(u,y_{1})K^{H^{\prime}}(v,y_{1})\,dy_{1}=R^{H^{\prime}}(u,v)$
and will be used repeatedly in the sequel.

To use this relation, we first expand the product in the expression for the
contraction in (\ref{aieaieaie}), taking care to keep track of the
indicator functions. The resulting initial expression for $(A_{i}\otimes
_{1}A_{i})(y_{1},y_{2})$ contains four terms, which are all of the following
form:%
\begin{eqnarray*}
C_{a,b} &:=& d ( H ) ^{2}\int_{0}^{1}dx\,1_{[0,a]} ( y_{1}\vee
x ) 1_{[0,b]} ( y_{2}\vee x ) \\
&&\hspace*{28.3pt}{} \times\int_{u=y_{1}\vee x}^{a}\partial_{1}K^{H^{\prime}} ( u,x )
\,\partial_{1}K^{H^{\prime}} ( u,y_{1} ) \,du\\
&&\hspace*{28.3pt}{} \times\int_{v=y_{2}\vee
x}^{b}\partial
_{1}K^{H^{\prime}} ( v,x )\,
\partial_{1}K^{H^{\prime}} ( v,y_{2} ) \,dv.
\end{eqnarray*}
Here, to perform a Fubini argument by bringing the integral over $x$ inside, we first
note that $x<u\wedge v$ while $u\in[ y_{1},a]$ and $v\in[
y_{2},b]$. Also, note that the conditions $x\leq u$ and $u\leq a$ imply
that $x\leq
a$ and thus $1_{[0,a]} ( y_{1}\vee x ) $ can be replaced, after
Fubini, by $1_{[0,a]} ( y_{1} ) $. Therefore, using (\ref{souvent}%
), the above expression equals%
\begin{eqnarray*}
C_{a,b} &=& d ( H ) ^{2}1_{[0,a]\times[0,b]} (
y_{1},y_{2} ) \int_{y_{1}}^{a}\partial_{1}K^{H^{\prime}} (
u,y_{1} ) \,du \int_{y_{2}}^{b}\partial_{1}K^{H^{\prime}} (
v,y_{2} ) \,dv\\
&&\hspace*{108.3pt}{} \times\int_{0}^{u\wedge v}\partial_{1}K^{H^{\prime}} (
u,x ) \,\partial_{1}K^{H^{\prime}} ( v,x ) \,dx\\
&=& d ( H ) ^{2}a(H)1_{[0,a]\times[0,b]} ( y_{1}%
,y_{2} ) \int_{u=y_{1}}^{a}\int_{v=y_{2}}^{b}\partial_{1}K^{H^{\prime}%
} ( u,y_{1} )\, \partial_{1}K^{H^{\prime}} ( v,y_{2} )\\
&&\hspace*{185.7pt}{}\times
\vert u-v \vert^{2H^{\prime}-2}\,du\,dv\\
&=& d ( H )
^{2}a(H)\int_{u=y_{1}}^{a}\int_{v=y_{2}}^{b}\partial_{1}K (
u,y_{1} ) \,\partial_{1}K^{H^{\prime}} ( v,y_{2} )
\vert u-v \vert^{2H^{\prime}-2}\,du\,dv.
\end{eqnarray*}
The last equality comes from the fact that the indicator
functions in $y_{1},y_{2}$ are redundant: they can be pulled back
into the integral over $du\,dv$ and, therein, the functions
$\partial_{1}K^{H^{\prime}} ( u,y_{1} ) $ and
$\partial_{1}K^{H^{\prime}} ( v,y_{2} ) $ are, by
definition as functions of $y_{1}$ and~$y_{2}$, supported by
smaller intervals than $[0,a]$ and $[0,b]$, namely $[0,u]$ and
$[0,v]$, respectively.

Now, the contraction $(A_{i}\otimes_{1}A_{i})(y_{1},y_{2})$ equals
$C_{i/N,i/N}+C_{ ( i-1 ) /N, ( i-1 ) /N}-C_{ (
i-1 ) /N,i/N}-C_{i/N, ( i-1 ) /N}$. Therefore, from the last
expression above,%
%
\begin{eqnarray}\label{contrac1}\hspace*{15pt}
&&(A_{i}\otimes_{1}A_{i})(y_{1},y_{2}) \nonumber\\
&&\qquad =a(H)d(H)^{2} \biggl( \int_{y_{1}%
}^{{i}/{N}}du\int_{y_{2}}^{{i}/{N}}dv\,\partial_{1}K^{H^{\prime}%
}(u,y_{1})\nonumber\\
&&\hspace*{173pt}\hspace*{-17pt}\hspace*{-66.5pt}{}\times\partial_{1}K^{H^{\prime}}(v,y_{2})|u-v|^{2H^{\prime}-2}
\nonumber\\
&&\qquad\quad\hspace*{59.9pt}{} -\int_{y_{1}}^{{i}/{N}}du\int_{y_{2}}^{({i-1})/{N}}dv\,\partial
_{1}K^{H^{\prime}}(u,y_{1})\nonumber\\
&&\hspace*{201pt}\hspace*{-16.39pt}\hspace*{-82.7pt}{}\times\partial_{1}K^{H^{\prime}}(v,y_{2}%
)|u-v|^{2H^{\prime}-2}\\
&&\qquad\quad\hspace*{59.9pt}{} -\int_{y_{1}}^{({i-1})/{N}}du\int_{y_{2}}^{{i}/{N}}dv\,\partial
_{1}K^{H^{\prime}}(u,y_{1})\nonumber\\
&&\hspace*{201pt}\hspace*{-16.61pt}\hspace*{-82.7pt}{}\times\partial_{1}K^{H^{\prime}}(v,y_{2}%
)|u-v|^{2H^{\prime}-2}\nonumber\\
&&\qquad\quad\hspace*{59.9pt}{} + \int_{y_{1}}^{({i-1})/{N}}du\int_{y_{2}}^{({i-1})/{N}%
}dv\,\partial_{1}K^{H^{\prime}}(u,y_{1})\nonumber\\
&&\hspace*{201pt}\hspace*{-51.1pt}{}\times\partial_{1}K^{H^{\prime}}%
(v,y_{2})|u-v|^{2H^{\prime}-2} \biggr) .\nonumber
\end{eqnarray}
Since the integrands in the above four integrals are identical, we can simplify
the above formula, grouping the first two terms, for instance, to
obtain an
integral of $v$ over $I_{i}= ( \frac{i-1}{N},\frac{i}{N} ] $, with
integration over $u$ in $[y_{1},\frac{i}{n}]$. The same operation on
the last
two terms gives the negative of the same integral over $v$, with
integration over $u$
in $[y_{1},\frac{i-1}{n}]$. Then, grouping these two resulting terms
yields a
single term, which is an integral for $ ( u,v ) $ over $I_{i}\times
I_{i}$. We obtain the following, final, expression for our contraction:%
%
\begin{eqnarray}\label{contrac11}\quad\quad
&&(A_{i}\otimes_{1}A_{i})(y_{1},y_{2})\nonumber\\[-8pt]\\[-8pt]
&&\qquad=a(H)d(H)^{2}
\iint_{I_{i}\times
I_{i}%
}\partial_{1}K^{H^{\prime}}(u,y_{1})\,\partial_{1}K^{H^{\prime}}(v,y_{2}%
)|u-v|^{2H^{\prime}-2}\,du\,dv.\nonumber
\end{eqnarray}

Now, since the integrands in the double Wiener integrals defining
$T_{2}$ are
symmetric, we get%
\[
\mathbf{E} [ T_{2}^{2} ] =N^{4H-2}16\cdot2!\sum_{i,j=1}^{N}\langle
A_{i}\otimes_{1}A_{i},A_{j}\otimes_{1}A_{j}\rangle_{L^{2}([0,1]^{2})}.
\]
To evaluate the inner product of the two contractions, we first use
Fubini with expression (\ref{contrac11}); by doing so, one must
realize that the support of $\partial_{1}K^{H^{\prime}}(u,y_{1})$ is
$ \{ u>y_{1} \} $, which then makes the upper limit $1$
for the integration in $y_{1}$ redundant; similar remarks hold with
respect to
$u^{\prime},v,v^{\prime}$ and $y_{2}$.
In other words, we have%
%
\begin{eqnarray}\label{aiaiai}\qquad
&& \langle A_{i}\otimes_{1}A_{i},A_{j}\otimes_{1}A_{j}\rangle_{L^{2}%
([0,1])^{2}}\nonumber\\
&&\qquad =a(H)^{2}d(H)^{4}\int_{0}^{1}\int_{0}^{1}dy_{1}\,dy_{2}\int_{I_{i}}%
\int_{I_{i}}\int_{I_{j}}\int_{I_{j}}du^{\prime}\,dv^{\prime}\,du\,dv\nonumber\\
&&\qquad\quad{}\times |u-v|^{2H^{\prime}-2}|u^{\prime}-v^{\prime}|^{2H^{\prime}-2}\,\partial
_{1}K^{H^{\prime}}(u,y_{1})\,\partial_{1}K^{H^{\prime}}(v,y_{2})\nonumber\\
&&\qquad\quad{}\times\partial
_{1}K^{H^{\prime}}(u^{\prime},y_{1})\,\partial_{1}K^{H^{\prime}}(v^{\prime
},y_{2})\nonumber\\
&&\qquad =a(H)^{2}d(H)^{4}\int_{I_{i}}\int_{I_{i}}\int_{I_{j}}\int_{I_{j}%
}|u-v|^{2H^{\prime}-2}|u^{\prime}-v^{\prime}|^{2H^{\prime}-2}\,du^{\prime
}\,dv^{\prime}\,dv\,du\\
&&\qquad\quad\hspace*{56.6pt}{}\times \int_{0}^{u\wedge u^{\prime}}\partial_{1}K^{H^{\prime}}(u,y_{1}%
)\,\partial_{1}K^{H^{\prime}}(u^{\prime},y_{1})\,dy_{1}\nonumber\\
&&\qquad\quad\hspace*{56.6pt}{}\times \int_{0}^{v\wedge
v^{\prime}}\partial_{1}K^{H^{\prime}}(v,y_{2})\,\partial_{1}K^{H^{\prime}%
}(v^{\prime},y_{2})\,dy_{2}\nonumber\\
&&\qquad =a(H)^{4}d(H)^{4}\int_{I_{i}}\int_{I_{i}}\int_{I_{j}}\int_{I_{j}%
}|u-v|^{2H^{\prime}-2}|u^{\prime}-v^{\prime}|^{2H^{\prime
}-2}|u-u^{\prime
}|^{2H^{\prime}-2}\nonumber\\
&&\hspace*{147pt}{}\times|v-v^{\prime}|^{2H^{\prime}-2}\,du^{\prime}\,dv^{\prime}\,dv\,du,\nonumber
\end{eqnarray}
where we have used the expression (\ref{souvent}) in the last step.
Therefore, we
immediately have %
%
\begin{eqnarray}\label{t2-1}\qquad
\mathbf{E} [ T_{2}^{2} ]
&=& N^{4H-2}32a(H)^{4}d(H)^{4}\nonumber\\
&&{}\times\sum
_{i,j=1}^{N}\int_{I_{i}}\int_{I_{i}}\int_{I_{j}}\int_{I_{j}}du^{\prime
}\,dv^{\prime}\,dv\,du\\
&&{} \times|u-v|^{2H^{\prime}-2}|u^{\prime}-v^{\prime}|^{2H^{\prime}%
-2}|u-u^{\prime}|^{2H^{\prime}-2}|v-v'|^{2H^{\prime
}-2}.\nonumber
\end{eqnarray}
By Lemma \ref{lemmaT2} in the \hyperref[Append]{Appendix}, we conclude that
%
\begin{eqnarray}\label{c3H}%
\lim_{N\rightarrow\infty}\mathbf{E} [ T_{2}^{2} ] N^{2-2H}%
&=& 64a(H)^{2}d(H)^{4} \biggl( \frac{1}{2H-1}-\frac{1}{2H} \biggr)
\nonumber\\
&=& 16d(H) ^{2}\\
&:=& c_{3,H}.\nonumber
\end{eqnarray}

\subsubsection{The term $T_{4}$}

Now, for the $L^{2}$-norm of the term denoted by
\[
T_{4}:=N^{2H-1}\sum_{i=1}^{N}I_{4}(A_{i}\otimes A_{i}),
\]
by the isometry formula for multiple stochastic integrals, and using a
correction term to account for the fact that the integrand in $T_{4}$ is
nonsymmetric, we have
\begin{eqnarray*}
\mathbf{E}[T_{4}^{2}] &=& 8N^{4H-2}\sum_{i,j=1}^{N}\langle A_{i}\otimes
A_{i};A_{j}\otimes A_{j}\rangle_{L^{2}([0,1]^{4})}\\
&&{} + 4N^{4H-2}\sum_{i,j=1}^{N}4\langle A_{i}\otimes_{1}A_{j};A_{j}\otimes
_{1}A_{i}\rangle_{L^{2}([0,1]^{2})}=:\mathcal{T}_{4,0}+\mathcal{T}_{4,1}.
\end{eqnarray*}
We separate the calculation of the two terms $\mathcal{T}_{4,0}$ and
$\mathcal{T}_{4,1}$ above. We will see that these two terms are exactly
of the
same magnitude, so both calculations must be performed precisely.

The first term, $\mathcal{T}_{4,0}$, can be written as%
\[
\mathcal{T}_{4,0}=8N^{4H-2}\sum_{i,j=1}^{N} \bigl\vert\langle A_{i}%
,A_{j}\rangle_{L^{2}([0,1]^{2})} \bigr\vert^{2}.
\]
We calculate each individual scalar product $\langle A_{i},A_{j}\rangle
_{L^{2}([0,1]^{2})}$ as%
\begin{eqnarray*}
&&\langle A_{i},A_{j}\rangle_{L^{2}([0,1]^{2})}\\
&&\qquad= \int_{0}^{1}\int_{0}^{1}
A_{i}(y_{1},y_{2})A_{j}(y_{1},y_{2})\,dy_{1}\,dy_{2} \\
&&\qquad= d(H)^{2}\int
_{0}^{1}\int
_{0}^{1}dy_{1}\,dy_{2}\,1_{[0,{i}/{N}\wedge{j}/{N}]}(y_{1}\vee
y_{2})\\
&&\qquad\quad\hspace*{27.8pt}{}\times \biggl( \int_{y_{1}\vee y_{2}}^{{i}/{N}}\partial_{1}K^{H^{\prime}%
}(u,y_{1})\,\partial_{1}K^{H^{\prime}}(u,y_{2})\,du-1_{[0,({i-1})/{N}]}%
(y_{1}\vee y_{2})\\
&&\qquad\quad\hspace*{27.8pt}\hspace*{96pt}{}\times\int_{y_{1}\vee y_{2}}^{({i-1})/{N}}\partial_{1}%
K^{H^{\prime}}(u,y_{1})\,\partial_{1}K^{H^{\prime}}(u,y_{2})\,du \biggr) \\
&&\qquad\quad\hspace*{27.8pt}{}\times \biggl( \int_{y_{1}\vee y_{2}}^{{j}/{N}}\partial_{1}K^{H^{\prime}%
}(v,y_{1})\partial_{1}K^{H^{\prime}}(v,y_{2})\,dv-1_{[0,({j-1})/{N}]}%
(y_{1}\vee y_{2})\\
&&\qquad\quad\hspace*{23.5pt}\hspace*{96pt}{}\times\int_{y_{1}\vee y_{2}}^{({j-1})/{N}}\,\partial_{1}%
K^{H^{\prime}}(v,y_{1})\,\partial_{1}K^{H^{\prime}}(v,y_{2})\,dv \biggr) \\
&&\qquad= d(H)^{2}\int_{({i-1})/{N}}^{{i}/{N}}\int_{({j-1})/{N}}^{{j}/{N}}du\,dv\,
\biggl[ \int_{0}^{u\wedge v}\partial_{1}K^{H^{\prime}}%
(u,y_{1})\,\partial_{1}K^{H^{\prime}}(v,y_{1})\,dy_{1} \biggr] ^{2}.
\end{eqnarray*}
Here, (\ref{souvent}) yields%
\[
\langle A_{i},A_{j}\rangle_{L^{2}([0,1]^{2})}=d(H)^{2}a(H)^{2}\int
_{I_{i}}%
\int_{I_{j}}|u-v|^{2H-2}\,du\,dv,
\]
where, we have again used the notation $I_{i}= ( \frac{i-1}{N},\frac{i}%
{N} ] $ for $i=1,\ldots,N$. We finally obtain
%
\begin{eqnarray}\label{prod1}\hspace*{25pt}
&&\langle A_{i},A_{j}\rangle_{L^{2}([0,1]^{2})} \nonumber\\[-8pt]\\[-8pt]
&&\qquad= \frac{d(H)^{2}a(H)^{2}}%
{H(2H-1)} \frac{1}{2}\biggl[ 2 \biggl\vert\frac{i-j}{N} \biggr\vert
^{2H}
- \biggl\vert\frac{i-j+1}{N} \biggr\vert^{2H}- \biggl\vert\frac{i-j-1}%
{N} \biggr\vert^{2H} \biggr],\nonumber
\end{eqnarray}
where, more precisely, $d(H)^{2}a(H)^{2}(H(2H-1))^{-1}=2$.
Specifically, with
the constants $c_{1,H}$, $c_{2,H}$ and $c_{1,H}^{\prime}$ given by
%
\begin{eqnarray}\label{c1H}%
c_{1,H} &:=& 2+\sum_{k=1}^{\infty} \bigl( 2k^{2H}- ( k-1 )
^{2H}- ( k+1 ) ^{2H} \bigr) ^{2},\nonumber\\
c_{2,H} &:=& 2H^{2} ( 2H-1 ) / ( 4H-3 ) ,\\
c_{1,H}^{\prime} &:=& \bigl(2H(2H-1)\bigr)^{2}=9/16,\nonumber
\end{eqnarray}
using Lemmas \ref{lemmac1}, \ref{lemmac2} and an analogous result for
$H=3/4$, we get, asymptotically for large $N$,%
%
\begin{eqnarray}\label{t1-1}%
\lim_{N\rightarrow\infty}N\mathcal{T}_{4,0} &=& 16c_{1,H},\qquad 1/2<H<\tfrac
{3}{4},
\\
\label{t1-2}%
\lim_{N\rightarrow\infty}N^{4-4H}\mathcal{T}_{4,0} &=& 16c_{2,H},\qquad
H>\tfrac{3}{4},
\\
\label{t1-3}%
\lim_{N\rightarrow\infty}\frac{N}{\log N}\mathcal
{T}_{4,0} &=& 16c_{1,H}^{\prime
}=16,\qquad
H=\frac{3}{4}.
\end{eqnarray}

The second term, $\mathcal{T}_{4,1}$, can be dealt with by obtaining an
expression for
\[
\langle A_{i}\otimes_{1}A_{j};A_{j}\otimes_{1}A_{i}\rangle
_{L^{2}([0,1]^{2})}%
\]
in the same way as the expression obtained in (\ref{aiaiai}). We get%
\begin{eqnarray*}
\mathcal{T}_{4,1} &=& 16N^{4H-2}\sum_{i,j=1}^{N}\langle A_{i}\otimes_{1}%
A_{j};A_{j}\otimes_{1}A_{i}\rangle_{L^{2}([0,1]^{2})}\\
&=& 16d(H)^{4}a(H)^{4}N^{-2}\sum_{i,j=1}^{N}\int_{0}^{1}\int_{0}^{1}\int
_{0}^{1}\int_{0}^{1}dy\,dz\,dy^{\prime}\,dz^{\prime}\\
&&\hspace*{110.1pt}{}\times |y-z+i-j|^{2H^{\prime}-2}|y^{\prime}-z^{\prime}+i-j|^{2H^{\prime}%
-2}\\
&&\hspace*{110.1pt}{}\times |y-y^{\prime}+i-j|^{2H^{\prime}-2}|z-z^{\prime}+i-j|^{2H^{\prime}-2}.
\end{eqnarray*}
Now, similarly to the proof of Lemma \ref{lemmaT2}, we find the the following
three asymptotic behaviors:
\begin{itemize}
\item if $H\in(\frac{1}{2},\frac{3}{4})$, then $\tau_{1,H}^{-1}%
N\mathcal{T}_{4,1}$ converges to $1$, where%
%
\begin{equation}\label{tau1}%
\tau_{1,H}:=16d(H)^{4}a(H)^{4}c_{1,H} ;
\end{equation}
\item if $H>\frac{3}{4}$, then $\tau_{2,H}^{-1}N^{4-4H}\mathcal{T}_{4,1}$
converges to $1$, where%
%
\begin{equation}\label{tau2}%
\tau_{2,H}:=32d(H)^{4}a(H)^{4}\int_{0}^{1}(1-x)x^{4H-4}\,dx;
\end{equation}
\item if $H=\frac{3}{4}$, then $\tau_{3,H}^{-1}(N/\log N)\mathcal{T}_{4,1}$
converges to $1$, where%
%
\begin{equation} \label{tau3}%
\tau_{3,H}:=32d(H)^{4}a(H)^{4}.
\end{equation}

\end{itemize}

Combining these results for $\mathcal{T}_{4,1}$ with those for $\mathcal
{T}%
_{4,0}$ in lines (\ref{t1-1}), (\ref{t1-2}) and (\ref{t1-3}), we obtain the
asymptotics of $\mathbf{E} [ T_{4}^{2} ] $ as $N\rightarrow\infty$:%
\begin{eqnarray*}
\lim_{N\rightarrow\infty}N\mathbf{E} [ T_{4}^{2} ]
&=&e_{1,H},\qquad
\mbox{if }H\in\bigl(\tfrac{1}{2},\tfrac{3}{4}\bigr);\\
\lim_{N\rightarrow\infty}N^{4-4H}\mathbf{E} [ T_{4}^{2} ]
&=&e_{2,H},\qquad
\mbox{if }H\in\bigl(\tfrac{3}{4},1\bigr);\\
\lim_{N\rightarrow\infty}\frac{N}{\log N}\mathbf{E} [ T_{4}%
^{2} ] &=& e_{3,H},\qquad
\mbox{if }H=\frac{3}{4},%
\end{eqnarray*}
where, with $\tau_{i,H}$, $i=1,2,3$, given in (\ref{tau1}), (\ref
{tau2}) and
(\ref{tau3}), we defined%
%
\begin{eqnarray}\label{e1H}%
e_{1,H}&:=&(1/2)c_{1,H}+\tau_{1,H},\nonumber\\
e_{2,H}&:=&(1/2)c_{2,H}+\tau
_{2,H},\\
e_{3,H}&:=&c_{3,H}+\tau_{3,H}.\nonumber
\end{eqnarray}

Taking into account the estimations (\ref{t1-1}), (\ref{t1-2}) and (\ref{t1-3}),
with $c_{3,H}$ in (\ref{c3H}), we see that $\mathbf{E} [ T_{4}%
^{2} ] $ is always of smaller order than $\mathbf{E} [ T_{2}%
^{2} ] $; therefore, the mean-square behavior of $V_{N}$ is given by that
of the term $T_{2}$ only, which means that we obtain, for every $H>1/2$,
%
\begin{equation}\label{v2}%
\lim_{N\rightarrow\infty}\mathbf{E} \biggl[ \biggl( N^{1-H}V_{N}(2,a)\frac
{1}{\sqrt{c_{3,H}}} \biggr) ^{2} \biggr] =1.
\end{equation}

\subsection{Normality of the fourth chaos term $T_{4}$ when $H\leq3/4$}

The calculations for $T_{4}$ above prove that $\lim_{N\rightarrow\infty
}\mathbf{E}[G_{N}^{2}]=1$ for $H<3/4$, where $e_{1,H}$ is given in (\ref{e1H})
and%
%
\begin{equation}\label{gN}%
G_{N}:=\sqrt{N}N^{2H-1}e_{1,H}^{-1/2}I_{4} \Biggl( \sum_{i=1}^{N}A_{i}\otimes
A_{i} \Biggr).\vadjust{\goodbreak}
\end{equation}
Similarly, for $H=\frac{3}{4}$, we showed that $\lim_{N\rightarrow\infty
}\mathbf{E}[\tilde{G}_{N}^{2}]=1$, where $e_{3,H}$ is given in (\ref
{e1H}) and
%
\begin{equation}\label{gNtilde}%
\tilde{G}_{N}:=\sqrt{\frac{N}{\log N}}N^{2H-1}e_{3,H}^{-1}I_{4} \Biggl(
\sum_{i=1}^{N}A_{i}\otimes A_{i} \Biggr) .
\end{equation}
Using the criterion of Nualart and Ortiz-Latorre [part (iv) in Theorem
\ref{NOTcrit}], we prove the following asymptotic normality for $G_{N}$ and
$\tilde{G}_{N}$.
\begin{theorem}
\label{geesbend}If $H\in(1/2,3/4)$, then $G_{N}$ given by (\ref{gN})
converges in distribution as%
%
\begin{equation} \label{convI4}%
\lim_{N\rightarrow\infty}G_{N}={\mathcal{N}}(0,1).
\end{equation}
If $H=3/4$, then $\tilde{G}_{N}$ given by (\ref{gNtilde}) converges in
distribution as%
%
\begin{equation}\label{convI5}%
\lim_{N\rightarrow\infty}\tilde{G}_{N}={\mathcal{N}}(0,1).
\end{equation}
\end{theorem}
\begin{pf}
We will denote by $c$ a generic positive constant not depending on $N$.
\setcounter{step}{-1}
\begin{step}[(Setup and expectation evaluation)]\label{step0}
Using the derivation
rule for multiple stochastic integrals, the Malliavin derivative of
$G_{N}$ is%
\[
D_{r}G_{N}=\sqrt{N}N^{2H-1}e_{1,H}^{-1/2}4\sum_{i=1}^{N}I_{3} \bigl(
(A_{i}\otimes A_{i})(\cdot,r) \bigr)
\]
and its norm is%
\begin{eqnarray*}
&&\Vert DG_{N}\Vert_{L^{2}([0,1])}^{2}\\
&&\qquad=N^{4H-1}16e_{1,H}^{-1}\sum_{i,j=1}%
^{N}\int_{0}^{1}dr\,I_{3} \bigl( (A_{i}\otimes A_{i})(\cdot,r) \bigr)
I_{3} \bigl( (A_{j}\otimes A_{j})(\cdot,r) \bigr) .
\end{eqnarray*}
The product formula (\ref{product}) gives
\begin{eqnarray*}
&&\Vert DG_{N}\Vert_{L^{2}([0,1])}^{2} \\
&&\qquad= N^{4H-1}16e_{1,H}^{-1}\\
&&\qquad\quad\hspace*{0pt}{}\times\sum
_{i,j=1}^{N}\int_{0}^{1}dr \,\bigl[
I_{6} \bigl( (A_{i}\otimes A_{i})(\cdot,r)\otimes(A_{j}\otimes A_{j}%
)(\cdot,r) \bigr) \\
&&\qquad\quad\hspace*{66.8pt}{} +9I_{4} \bigl( (A_{i}\otimes A_{i})(\cdot,r)\otimes_{1}(A_{j}\otimes
A_{j})(\cdot,r) \bigr) \\
&&\qquad\quad\hspace*{66.8pt}{} +9I_{2} \bigl( (A_{i}\otimes A_{i})(\cdot,r)\otimes_{2}(A_{j}\otimes
A_{j})(\cdot,r) \bigr) \\
&&\qquad\quad\hspace*{66.8pt}{} +3!I_{0} \bigl( (A_{i}\otimes A_{i})(\cdot,r)\otimes_{3}%
(A_{j}\otimes A_{j})(\cdot,r) \bigr)
\bigr]\\
&&\qquad=:J_{6}+J_{4}+J_{2}+J_{0}.
\end{eqnarray*}
First, note that, for the nonrandom term $J_{0}$ that gives the
expected value
of the above, we have%
\begin{eqnarray*}
J_{0} &=& 16e_{1,H}^{-1}N^{4H-1}3!\sum_{i,j=1}^{N}\int_{[0,1]^{4}}A_{i}%
(y_{1},y_{2})A_{i}(y_{3},y_{4})A_{j}(y_{1},y_{2})\\
&&\hspace*{119.5pt}{}\times A_{j}(y_{3},y_{4}%
)\,dy_{1}\,dy_{2}\,dy_{3}\,dy_{4}\\
&=& 96N^{4H-1}e_{1,H}^{-1}\sum_{i,j=1}^{N} \bigl\vert\langle A_{i}%
,A_{j}\rangle_{_{L^{2}([0,1]^{2})}} \bigr\vert^{2}.
\end{eqnarray*}
This sum has already been treated: we know from (\ref{t1-1}) that $J_{0}/4$
converges to~$1$, that is, that $\lim_{_{N\rightarrow\infty}}
\mathbf{E}[\Vert DG_{N}\Vert_{L^{2}([0,1])}^{2}]=4$. This means, by the
Nualart--Ortiz-Latorre criterion, that we only need to show that all other
terms $J_{6},J_{4},J_{2}$ converge to zero in $L^{2}(\Omega)$ as
$N\rightarrow\infty$.
\end{step}
\begin{step}[(Order-6 chaos term)]\label{step1}
We first consider the
term $J_{6}%
$:%
\begin{eqnarray*}
J_{6}&=&
cN^{4H-1}\sum_{i,j=1}^{N}\int_{0}^{1}dr\,I_{6} \bigl( (A_{i}\otimes
A_{i})(\cdot,r)\otimes\bigl(A_{j}\otimes A_{j}(\cdot,r)\bigr) \bigr) \\
&=&
cN^{4H-1}%
\sum_{i,j=1}^{N}I_{6} \bigl( (A_{i}\otimes A_{j})\otimes(A_{i}\otimes_{1}%
A_{j}) \bigr) .
\end{eqnarray*}
We study the mean square of this term. We have, since the $L^{2}$-norm
of the
symmetrization is less than the $L^{2}$-norm of the corresponding
unsymmetrized function,
\begin{eqnarray*}
&& \mathbf{E} \Biggl[ \Biggl( \sum_{i,j=1}^{N}I_{6} \bigl( (A_{i}\otimes
A_{j})\otimes(A_{i}\otimes_{1}A_{j}) \bigr) \Biggr) ^{2} \Biggr] \\
&&\qquad \leq6!\sum_{i,j,k,l}\langle(A_{i}\otimes A_{j})\otimes(A_{i}\otimes
_{1}A_{j}),(A_{k}\otimes A_{l})\otimes(A_{k}\otimes_{1}A_{l})\rangle
_{L^{2}([0,1]^{6})}\\
&&\qquad =6!\sum_{i,j,k,l}\langle A_{i},A_{k}\rangle_{L^{2}([0,1]^{2})}\langle
A_{j},A_{l}\rangle_{L^{2}([0,1]^{2})}\langle A_{i}\otimes_{1}A_{j}%
,A_{k}\otimes_{1}A_{l}\rangle_{L^{2}([0,1]^{2})}.
\end{eqnarray*}
We get
\begin{eqnarray*}
\mathbf{E} [ J_{6}^{2} ]
&\leq& cN^{8H-2}\sum_{i,j,k,l}\int_{I_{i}}\,du\int_{I_{j}}dv\int_{I_{k}%
}du^{\prime}\int_{I_{l}}dv^{\prime}\\
&&\hspace*{59.2pt}{}\times|u-v|^{2H^{\prime}-2}|u-u^{\prime
}|^{2H^{\prime}-2}|v-v^{\prime}|^{2H^{\prime}-2}|u^{\prime}-v^{\prime
}|^{2H^{\prime}-2}\\
&&\hspace*{59.2pt}{}\times\biggl[ 2 \biggl\vert\frac{i-k}{N} \biggr\vert^{2H}- \biggl\vert
\frac{i-k+1}{N} \biggr\vert^{2H}- \biggl\vert\frac{i-k-1}{N} \biggr\vert
^{2H} \biggr]\\
&&\hspace*{59.2pt}{}\times\biggl[ 2 \biggl\vert\frac{j-l}{N} \biggr\vert^{2H}- \biggl\vert
\frac{j-l+1}{N} \biggr\vert^{2H}- \biggl\vert\frac{j-l-1}{N} \biggr\vert
^{2H} \biggr] .
\end{eqnarray*}
First, we show that for $H\in(1/2,3/4)$, we have, for large $N$,
%
\begin{equation}\label{EJ1}%
\mathbf{E} [ J_{6}^{2} ] \leq cN^{8H-6}.
\end{equation}

With the notation as in Step \ref{step1} of this proof, making the change of variables
$\bar{u}=(u-\frac{i-1}{N})N$, and similarly for the other integrands, we
obtain
\begin{eqnarray*}
\mathbf{E} [ J_{6}^{2} ] &\leq& cN^{8H-2}\frac{1}{N^{8H^{\prime
}-8}}\frac{1}{N^{4}}\frac{1}{N^{4H}}\\
&&{}\times\sum_{i,j,k,l}\int_{[0,1]^{4}%
}du\,dv\,du^{\prime}\,dv^{\prime}\\
&&\hspace*{33.47pt}{}
\times|u-v+i-j|^{2H^{\prime}-2}|u-u^{\prime}+i-k|^{2H^{\prime}-2}\\
&&\hspace*{33.47pt}{}\times
|u^{\prime}-v+j-k|^{2H^{\prime}-2}|v-v^{\prime}+k-l|^{2H^{\prime}-2}\\
&&\hspace*{33.47pt}{} \times( 2 \vert i-k \vert^{2H}- \vert i-k+1 \vert
^{2H}- \vert i-k-1 \vert^{2H} )\\
&&\hspace*{33.47pt}{}\times ( 2 \vert
j-l \vert^{2H}- \vert j-l+1 \vert^{2H}- \vert
j-l-1 \vert^{2H} ) \\
&=& c\frac{1}{N^{2}}\sum_{i,j,k,l}\int_{[0,1]^{4}}du\,dv\,du^{\prime
}\,dv^{\prime}\\
&&\hspace*{43.47pt}{}
\times|u-v+i-j|^{2H^{\prime}-2}|u-u^{\prime}+i-k|^{2H^{\prime}-2}\\
&&\hspace*{43.47pt}{}\times
|u^{\prime}-v+j-k|^{2H^{\prime}-2}|v-v^{\prime}+k-l|^{2H^{\prime}-2}\\
&&\hspace*{43.47pt}{} \times( 2 \vert i-k \vert^{2H}- \vert i-k+1 \vert
^{2H}- \vert i-k-1 \vert^{2H} )\\
&&\hspace*{43.47pt}{}\times( 2 \vert
j-l \vert^{2H}- \vert j-l+1 \vert^{2H}- \vert
j-l-1 \vert^{2H} ).
\end{eqnarray*}
Again, we use the fact that the dominant part in the above expression
is the
one in where all indices are distant by at least two units. In this
case, up to a
constant, we have the upper bound $|i-k|^{2H-2}$ for the quantity $ (
2 \vert i-k \vert^{2H}- \vert i-k+1 \vert^{2H}- \vert
i-k-1 \vert^{2H} ) $. By using Riemann sums, we can write
\[
\mathbf{E} [ J_{6}^{2} ] \leq c\frac{1}{N^{2}}N^{4} \biggl( \frac
{1}{N^{4}}\sum_{i,j,k,l}f\biggl(\frac{i}{N},\frac{j}{N},\frac{k}{N},\frac{l}%
{N}\biggr) \biggr) N^{8H^{\prime}-8}N^{4H-4},%
\]
where $f$ is a Riemann integrable function on $[0,1]^{4}$ and the
Riemann sum
converges to the finite integral of $f$ therein. Estimate (\ref{EJ1})
follows.
\end{step}
\begin{step}[(Chaos terms of orders 4 and 2)]\label{step2}
To treat the term
\[
J_{4}=cN^{4H-1}\sum_{i,j=1}^{N}\int_{0}^{1}dr\,I_{4} \bigl( (A_{i}\otimes
A_{i})(\cdot,r)\otimes_{1}(A_{j}\otimes A_{j})(\cdot,r) \bigr) ,
\]
since $I_{4}(g)=I_{4}(\tilde{g})$, where $\tilde{g}$ denotes the symmetrization
of the function $g$, we can write
\begin{eqnarray*}
J_{4} &=& cN^{4H-1}\sum_{i,j=1}^{N}\langle A_{i},A_{j}\rangle_{L^{2}%
(0,1]^{2}}I_{4} ( A_{i}\otimes A_{j} )\\
&&{} +cN^{4H-1}I_{4}\sum
_{i,j=1}^{N}(A_{i}\otimes_{1}A_{j})\otimes(A_{i}\otimes_{1}A_{j}%
)\\
&=:&J_{4,1}+J_{4,2}.
\end{eqnarray*}
Both terms above have been treated in previous computations. To
illustrate it,
the first summand $J_{4,1}$ can be bounded above as follows:
\begin{eqnarray*}
\mathbf{E} \vert J_{4,1} \vert^{2}
&\leq& cN^{8H-2}\sum_{i,j,k,l=1}^{N}\langle A_{i},A_{j}\rangle_{L^{2}%
([0,1]^{2})}\langle A_{i},A_{k}\rangle_{L^{2}([0,1]^{2})}\\
&&\hspace*{69.2pt}{}\times\langle A_{k}%
,A_{l}\rangle_{L^{2}([0,1]^{2})}\langle A_{j},A_{l}\rangle
_{L^{2}([0,1]^{2}%
)}\\
&=& cN^{8H-2}\sum_{i,j,k,l=1}^{N} \biggl[ \biggl( \frac{i-j+1}{N} \biggr)
^{2H} + \biggl( \frac{i-j-1}{N} \biggr) ^{2H}-2 \biggl( \frac{i-j}{N} \biggr)
^{2H} \biggr] \\
&&\hspace*{70.1pt}{}\times \biggl[ \biggl( \frac{i-k+1}{N} \biggr) ^{2H}\\
&&\hspace*{88.1pt}{}+ \biggl( \frac{i-k-1}%
{N} \biggr) ^{2H}-2 \biggl( \frac{i-k}{N} \biggr) ^{2H} \biggr] \\
&&\hspace*{70.1pt}{}\times \biggl[ \biggl(
\frac{j-l+1}{N} \biggr) ^{2H}\\
&&\hspace*{88.1pt}{} + \biggl( \frac{j-l-1}{N} \biggr) ^{2H}-2 \biggl(
\frac{j-l}{N} \biggr) ^{2H} \biggr] \\
&&\hspace*{70.1pt}{}\times \biggl[ \biggl( \frac{k-l+1}{N} \biggr) ^{2H}+ \biggl( \frac{k-l-1}%
{N} \biggr) ^{2H}-2 \biggl( \frac{k-l}{N} \biggr) ^{2H} \biggr]
\end{eqnarray*}
and, using the same bound $c|i-j|^{2H-2}$ for the quantity $|i-j+1|^{2H}
+|i-j-1|^{2H}-2|i-j|^{2H}$ when $|i-j|\geq2$, we obtain
\begin{eqnarray*}
\mathbf{E} \vert J_{4,1} \vert^{2}
&\leq&
cN^{8H-2}N^{-8H}%
\sum_{i,j,k,l=1}^{N}|i-j|^{2H-2}|i-k|^{2H-2}|j-l|^{2H-2}|k-l|^{2H-2}\\
&\leq&
cN^{8H-6}\frac{1}{N^{4}}\sum_{i,j,k,l=1}^{2}\frac{|i-j|^{2H-2}%
|i-k|^{2H-2}|j-l|^{2H-2}|k-l|^{2H-2}}{N^{4(2H-2)}}.
\end{eqnarray*}
This tends to zero at the speed $N^{8H-6}$ as $N\rightarrow\infty$ by a
Riemann sum argument since $H<\frac{3}{4}$.

One can also show that $\mathbf{E} \vert J_{4,2} \vert^{2}$
converges to zero at the same speed because
\begin{eqnarray*}
\mathbf{E} \vert J_{4,2} \vert^{2}
&=&
cN^{8H-2}\sum_{i,j,k,l=1}%
^{N}\langle(A_{i}\otimes_{1}A_{j}),(A_{k}\otimes_{1}A_{l})\rangle
_{L^{2}([0,1]^{2})}^{2}\\
&\leq& N^{8H-2}N^{-2(8H^{\prime}-8)}N^{-8}\\
&&\hspace*{0pt}{}\times \sum_{i,j,k,l=1}^{N} \biggl( \int_{[0,1]^{4}} ( |u-v+i-j|\\
&&\hspace*{80.7pt}{}\times| u^{\prime
}-v^{\prime}+k -l|\\
&&\hspace*{80.7pt}{}\times| u -u^{\prime}+i-k|\\
&&\hspace*{83pt}{}\times| v-v^{\prime}+j-l|) ^{2H^{\prime}%
-2}\,dv'\,du^{\prime}\,dv\,du \biggr) ^{2}\\
&\leq& cN^{8H-6}.
\end{eqnarray*}
Thus, we obtain
%
\begin{equation}\label{EJ2}%
\mathbf{E} [ J_{4}^{2} ] \leq cN^{8H-6}.
\end{equation}

A similar behavior can be obtained for the last term $J_{2}$ by
repeating the
above arguments:
%
\begin{equation} \label{EJ3}%
\mathbf{E} [ J_{2}^{2} ] \leq c N^{8H-6}.
\end{equation}
\end{step}
\begin{step}[(Conclusion)]\label{step3}
Combining (\ref{EJ1}), (\ref{EJ2}) and
(\ref{EJ3}) and recalling the convergence result for $\mathbf{E}%
[ T_{4}^{2} ] $ proven in the previous subsection, we can apply
the Nualart--Ortiz-Latorre criterion and use the same method as in the case
$H<\frac{3}{4}$ for $H=3/4$, to conclude the proof.\qed
\end{step}
\noqed\end{pf}

\subsection{Nonnormality of the second chaos term $T_{2}$ and limit of the
$2$-variation}

This paragraph studies the asymptotic behavior of the term denoted by $T_{2}$
which appears in the decomposition of $V_{N}(2,a)$. Recall that this is the
dominant term, given by%
\[
T_{2}=4N^{2H-1}I_{2} \Biggl(
\sum_{i=1}^{N}A_{i}\otimes_{1}A_{i} \Biggr),
\]
and, with $\sqrt{c_{3,H}}=4d ( H ) $ given in (\ref{c3H}), we
have shown that
\[
\lim_{N\rightarrow\infty}\mathbf{E} [ ( N^{1-H}T_{2}c_{3,H}%
^{-1/2} ) ^{2} ] =1.
\]
With $T_{N}:=N^{1-H}T_{2}c_{3,H}^{-1/2}$, one can show that in
$L^{2}(\Omega)$,
\[
\lim_{N\rightarrow\infty}\Vert
DT_{N}\Vert_{L^{2}([0,1])}^{2}=2+c,
\]
where $c$ is a strictly positive
constant. As a consequence, the Nualart--Ortiz-Latorre criterion can be used
to deduce that the $T_{N}$ do not converge to the standard normal law.
However, it is straightforward to find the limit of $T_{2}$, and
thus of $V_{N}$, in $L^{2} ( \Omega)$, in this case. We
have the following result.
\begin{theorem}
\label{T2toRose} For all $H\in(1/2,1)$, the normalized $2$-variation
$N^{1-H}V_{N}(2,\break a)/ ( 4d ( H ) ) $ converges in
$L^{2} ( \Omega) $ to the Rosenblatt random variable $Z (
1 ) $. Note that this is the actual observed value of the Rosenblatt
process at time $1$.
\end{theorem}
\begin{pf}
Since we already proven that $N^{1-H}T_{4}$ converges to $0$ in $L^{2} (
\Omega) $, it is sufficient to prove that $N^{1-H}T_{2}/ (
4d ( H ) ) -Z ( 1 ) $ converges to $0$ in
$L^{2} ( \Omega) $. Since $T_{2}$ is a second-chaos random
variable, that is, is of the form $I_{2} ( f_{N} )$, where $f_{N}$ is a
symmetric function in $L^{2} ( [0,1]^{2} ) $, it is sufficient to
prove that
\[
\frac{N^{1-H}}{4d(H)}f_{N}%
\]
converges to $L_{1}$ in $L^{2} ( [0,1]^{2} ) $, where $L_{1}$ is
given by (\ref{lt}). From (\ref{contrac11}), we get
%
\begin{eqnarray}\label{fNtrue}\qquad
f_{N} ( y_{1},y_{2} ) &=& 4N^{2H-1}a(H)d(H)^{2}\nonumber\\
&&{}\times \sum_{i=1}%
^{N} \biggl( \iint_{I_{i}\times I_{i}}|u-v|^{2H^{\prime}-2}\\
&&\hspace*{65.99pt}{}\times\partial
_{1}K^{H^{\prime}}(u,y_{1})\,\partial_{1}K^{H^{\prime}}(v,y_{2})\,du\,dv
\biggr).\nonumber
\end{eqnarray}

We now show that $\frac{N^{1-H}}{4d(H)}f_{N}$ converges pointwise, for
$y_{1},y_{2}\in[0,1]$, to the kernel of the Rosenblatt random variable.
On the interval $I_{i}\times I_{i}$, we may replace the evaluation of
$\partial_{1}K^{H^{\prime}}$ and $\partial_{1}K^{H^{\prime}}$ at $u$
and $v$
by setting $u=v=i/N$. We then get that $f_{N} ( y_{1},y_{2} ) $ is
asymptotically equivalent to
\begin{eqnarray*}
&& 4N^{2H-1}a ( H ) d ( H ) ^{2}\sum_{i=1}^{N}%
\mathbf{1}_{i/N\geq y_{1}\vee y_{2}}\,\partial_{1}K^{H^{\prime}}(i/N,y_{1}
)\,\partial_{1}K^{H^{\prime}}(i/N,y_{2})\\
&&\hspace*{95.3pt}\quad{}\times\iint_{I_{i}\times
I_{i}}du\,dv\,
\vert
u-v \vert^{2H^{\prime}-2}\\
&&\qquad =4N^{H-1}d ( H ) ^{2}\frac{1}{N}\sum_{i=1}^{N}\mathbf{1}%
_{i/N\geq y_{1}\vee y_{2}}\,\partial_{1}K^{H^{\prime}}(i/N,y_{1})\,\partial
_{1}K^{H^{\prime}}(i/N,y_{2}),
\end{eqnarray*}
where we have used the identity $\iint_{I_{i}\times I_{i}}du\,dv\, \vert
u-v \vert^{2H^{\prime}-2}=a ( H ) ^{-1}N^{-2H^{\prime}%
}=a ( H ) ^{-1}N^{-H-1}$. Therefore, we can write, for every
$y_{1},y_{2}\in(0,1)^{2}$, by\vadjust{\goodbreak} invoking a Riemann sum approximation,%
\begin{eqnarray*}
&&\lim_{N\rightarrow\infty}\frac{N^{1-H}}{4d(H)}f_{N} ( y_{1},y_{2} )
\\
&&\qquad =d(H)\lim_{N\rightarrow\infty}\frac{1}{N}\sum_{i=1}^{N}\mathbf
{1}_{i/N\geq
y_{1}\vee y_{2}}\,\partial_{1}K^{H^{\prime}}(i/N,y_{1})\,\partial
_{1}K^{H^{\prime
}}(i/N,y_{2})\\
&&\qquad =d(H)\int_{y_{1}\vee y_{2}}^{1}\partial_{1}K^{H^{\prime}}(u,y_{1}%
)\partial_{1}K^{H^{\prime}}(u,y_{2})\,du=L_{1} ( y_{1,}y_{2} ).
\end{eqnarray*}

To complete the proof, it suffices to check that the sequence
$(4d(H))^{-1}\times\break N^{1-H} f_{N}$ is Cauchy in $L^{2}([0,1]^{2})$ [indeed,
this implies that $(4d(H))^{-1}N^{1-H}f_{N}$ has a limit in
$L^{2}([0,1]^{2})$, which obviously coincides with the a.e. limit
$L_{1}$ and then the multiple integral
$I_{2}((4d(H))^{-1}N^{1-H}f_{N})$ will converge to $I_{2}(L_{1})$].
This can be checked by means of a straightforward calculation. Indeed, one
has, with $C(H)$ a positive constant not depending on $M$ and $N$,
%
\begin{eqnarray}\label{forcauchy}
&&\hspace*{-5.9pt} \Vert N^{1-H}f_{N}-M^{1-H}f_{M}\Vert_{L^{2}([0,1]^{2})}^{2}\nonumber
\\[-0.7pt]
&&\hspace*{-8pt}\qquad =C(H)N^{2H}\sum_{i,j=1}^{N}\int_{I_{i}}\int_{I_{i}}\int_{I_{j}}\int
_{I_{j}%
}|u-v|^{2H^{\prime}-2}|u^{\prime}-v^{\prime}|^{2H^{\prime
}-2}\nonumber\\[-0.7pt]
&&\hspace*{-8pt}\hspace*{160.3pt}{}\times|u-u^{\prime
}|^{2H^{\prime}-2}|v-v^{\prime}|^{2H^{\prime}-2}\,du^{\prime}\,dv^{\prime
}\,du\,dv\nonumber\\[-0.7pt]
&&\hspace*{-8pt}\qquad\quad{} +C(H)M^{2H}\nonumber\\[-0.7pt]
&&\hspace*{-8pt}\hspace*{43.2pt}{}\times\hspace*{-1.5pt}\sum_{i,j=1}^{M}\hspace*{-1pt}\int_{({i-1})/{M}}^{{i}/{M}}\hspace*{-0.5pt}\int
_{({i-1})/{M}}^{{i}/{M}}\hspace*{-0.5pt}\int_{({j-1})/{M}}^{{j}/{M}}\hspace*{-0.5pt}\int
_{({j-1})/{M}}^{{j}/{M}}\hspace*{-0.5pt}|u-v|^{2H^{\prime}-2}|u^{\prime}-v^{\prime
}|^{2H^{\prime}-2}\nonumber\\[-0.7pt]
&&\hspace*{-8pt}\hspace*{231.5pt}{}\times\hspace*{-0.3pt}|u-u^{\prime}|^{2H^{\prime}-2}\\[-0.7pt]
&&\hspace*{-8pt}\hspace*{231.5pt}{}\times\hspace*{-0.3pt}|v-v^{\prime
}|^{2H^{\prime}%
-2}\,du^{\prime}\,dv^{\prime}\,du\,dv\nonumber\\[-0.7pt]
&&\hspace*{-8pt}\qquad\quad{} -2C(H)M^{1-H}N^{1-H}M^{2H-1}N^{2H-1}\nonumber\\[-0.7pt]
&&\hspace*{-8pt}\hspace*{43.2pt}{}\times\sum_{i=1}^{N}\sum_{j=1}^{M}\int
_{I_{i}}\int_{I_{i}}\int_{({j-1})/{M}}^{{j}/{M}}\int_{
({j-1})/{M}%
}^{{j}/{M}}du^{\prime}\,dv^{\prime}\,du\,dv\nonumber\\[-0.7pt]
&&\hspace*{-8pt}\hspace*{43.2pt}{} \times|u-v|^{2H^{\prime}-2}|u^{\prime}-v^{\prime}|^{2H^{\prime}
-2}\nonumber\\[-0.7pt]
&&\hspace*{-8pt}\hspace*{43.2pt}{}
\times|u-u^{\prime}|^{2H^{\prime}-2}|v-v^{\prime}|^{2H^{\prime}-2}.\nonumber
\end{eqnarray}
The first two terms have already been studied in Lemma \ref{lemmaT2}.
We have
shown that
\begin{eqnarray*}
&&
N^{2H}\sum_{i,j=1}^{N}\int_{I_{i}}\int_{I_{i}}\int_{I_{j}}\int_{I_{j}%
}|u-v|^{2H^{\prime}-2}|u^{\prime}-v^{\prime}|^{2H^{\prime
}-2}\\
&&\qquad\hspace*{78.4pt}{}\times |u-u^{\prime
}|^{2H^{\prime}-2}|v-v^{\prime}|^{2H^{\prime}-2}\,du^{\prime}\,dv^{\prime}\,du\,dv
\end{eqnarray*}
converges to $(a ( H ) ^{2}H(2H-1))^{-1}$. Thus, each of the first
two terms in (\ref{forcauchy}) converge to $C ( H ) $ times that
same constant as $M,N$ go to infinity. By the change of variables which
has already been used
several times, $\bar{u}=(u-\frac{i}{N})N$, the last term in (\ref{forcauchy})
is equal to%
\begin{eqnarray*}
&& C ( H ) (MN)^{H}\frac{1}{N^{2}M^{2}}(NM)^{2H^{\prime}-2}\\
&&\qquad{}\times
\sum_{i=1}^{N}\sum_{j=1}^{M}\int_{[0,1]^{4}}du\,dv\,du^{\prime}\,dv^{\prime}\\
&&\hspace*{92pt}\hspace*{-27.8pt}{} \times|u-v|^{2H^{\prime}-2}|u^{\prime}-v^{\prime}|^{2H^{\prime}%
-2}\\
&&\hspace*{92pt}\hspace*{-27.8pt}{} \times \biggl\vert\frac{u}{N}-\frac{u^{\prime}}{M}+\frac{i}{N}-\frac{j}%
{M} \biggr\vert^{2H^{\prime}-2}\\
&&\hspace*{92pt}\hspace*{-27.8pt}{} \times \biggl\vert\frac{v}{N}-\frac{v^{\prime}}%
{M}+\frac{i}{N}-\frac{j}{M} \biggr\vert^{2H^{\prime}-2}\\
&&\qquad =\frac{C ( H ) }{MN}\sum_{i=1}^{N}\sum_{j=1}^{M}\int_{[0,1]^{4}%
}du\,dv\,du^{\prime}\,dv^{\prime}\\
&&\hspace*{120pt}\hspace*{-28.2pt}{} \times|u-v|^{2H^{\prime}-2}|u^{\prime}-v^{\prime}|^{2H^{\prime}-2}
\\
&&\hspace*{120pt}\hspace*{-28.2pt}{}\times \biggl\vert\frac{u}{N}-\frac{u^{\prime}}{M}+\frac{i}{N}-\frac{j}%
{M} \biggr\vert^{2H^{\prime}-2}
\\
&&\hspace*{120pt}\hspace*{-28.2pt}{}\times \biggl\vert\frac{v}{N}-\frac{v^{\prime}}%
{M}+\frac{i}{N}-\frac{j}{M} \biggr\vert^{2H^{\prime}-2}.
\end{eqnarray*}
For large $i,j$, the term $\frac{u}{N}-\frac{u^{\prime}}{M}$ in
front of $\frac{i}{N}-\frac{j}{M}$ is negligible and can be ignored.
Therefore, the last
term in (\ref{forcauchy}) is equivalent to a Riemann sum than tends, as
$M,N\rightarrow\infty$, to the constant $ ( \int_{0}^{1}\int_{0}%
^{1}|u-v|^{2H^{\prime}-2}\,du\,dv ) ^{2}\int_{0}^{1}\int_{0}^{1}%
|x-y|^{2(2H^{\prime}-2)}$. This is precisely equal to $2(a ( H )
^{2}H(2H-1))^{-1}$, that is, the limit of the sum of the first two
terms in
(\ref{forcauchy}). Since the last term has a leading negative sign, the
announced Cauchy convergence is established, completing the proof of the
theorem.
\end{pf}
\begin{remark}
\label{rem2}One can show that the 2-variations $V_{N}(2,a)$ converge to
zero almost surely as $N$ goes to infinity. Indeed, the results in this
section already show that $V_{N} ( 2,a ) $ converges to $0$ in
$L^{2} ( \Omega) $, and thus in probability, as $N\rightarrow
\infty$; the almost sure convergence is obtained by using an argument in
\cite{coeur} (proof of Proposition 1) based on Theorem 6.2 in \cite{Doob}
which gives the equivalence between the almost sure convergence and the
mean-square convergence for empirical means of discrete stationary
processes. This almost-sure convergence can also be proven by hand in the
following standard way. Since $V_{N} ( 2,a ) $ is in the fourth
Wiener chaos, it is known that its $2p$th moment is bounded above by
$c_{p} ( \mathbf{E} [ ( V_{N} ( 2,a ) )
^{2} ] ) ^{p/2}$, where $c_{p}$ depends only on $p$. By choosing
$p$ large enough, via Chebyshev's inequality, the Borel--Cantelli lemma
yields the desired conclusion.
\end{remark}

\subsection{Normality of the adjusted variations}

According to Theorem \ref{T2toRose}, which we just proved, in the
Rosenblatt case, the standardization of the random variable
$V_{N}(2,a)$ does not converge to the normal law. But, this
statistic, which can be written as $V_{N}=T_{4}+T_{2}$, has a small
\textit{normal part}, which is given by the asymptotics of the term
$T_{4}$, as we can see from Theorem \ref{geesbend}. Therefore,
$V_{N}-T_{2}$ will converge (under suitable scaling) to the Gaussian
distribution. Of course, the term~$T_{2}$, which is an iterated
stochastic integral, is not practical because it cannot be observed.
But, replacing it with its limit $Z(1)$ (this is observed), one can
define an adjusted version of the statistic $V_{N}$ that
converges, after standardization, to the standard normal law.

The proof of this fact is somewhat delicate. If we are to subtract a multiple
of $Z ( 1 ) $ from $V_{N}$ in order to recuperate $T_{4}$ and hope
for a normal convergence, the first calculation would have to be as
follows:%
%
\begin{eqnarray}\label{VNwantnormal}
V_{N}(2,a)-\frac{\sqrt{c_{3,H}}}{N^{1-H}}Z(1) &=& V_{N}(2,a)-T_{2}+T_{2}
-\frac{\sqrt{c_{3,H}}}{N^{1-H}}Z(1)\nonumber\\
&=& T_{4}+\frac{\sqrt{c_{3,H}}}{N^{1-H}} \biggl[ \frac{N^{1-H}}{\sqrt{c_{3,H}}
}T_{2}-Z(1) \biggr] \\
&:=& T_{4}+U_{2}.\nonumber
\end{eqnarray}
The term $T_{4}$, when normalized as $\frac{\sqrt{N}}{\sqrt{e_{1,H}}}T_{4}$,
converges to the standard normal law, by Theorem \ref{geesbend}.
To get a normal convergence for the entire expression in (\ref{VNwantnormal}),
one may hope that the additional term $U_{2}:=\frac{\sqrt
{c_{3,H}}}{N^{1-H}%
} [ \frac{N^{1-H}}{\sqrt{c_{3,H}}}T_{2}-Z(1) ] $ goes to $0$
``fast enough.'' It is certainly true that
$U_{2}$ does go to $0$, as we just saw in Theorem \ref{T2toRose}.
However, the proof of that theorem did not investigate the speed of this
convergence of $U_{2}$. For this convergence to be ``fast
enough,'' one must multiply the expression by the rate
$\sqrt{N}$ which is needed to ensure the normal convergence of $T_{4}$: we
would need $U_{2}\ll N^{-1/2}$. Unfortunately, this is not true. A more
detailed calculation will show that $U_{2}$ is precisely of order $\sqrt{N}$.
This means that we should investigate whether $\sqrt{N}U_{2}$ itself converges
in distribution to a normal law. Unexpectedly, this turns out to be
true if
(and only if) $H<2/3$.
\begin{proposition}
\label{AjustPro}With $U_{2}$ as defined in (\ref{VNwantnormal}) and $H<2/3$,
we have that $\sqrt{N}U_{2}$ converges in distribution to a centered
normal with
variance equal to%
%
\begin{equation} \label{f1H}%
f_{1,H}:=32d ( H ) ^{4}a ( H ) ^{2}\sum_{k=1}^{\infty
}k^{2H-2}F \biggl( \frac{1}{k} \biggr),
\end{equation}
where the function $F$ is defined by
%
\begin{eqnarray} \label{eff}%
&& F(x)=\int_{[0,1]^{4}}du\,dv\,du^{\prime}\,dv^{\prime}|(u-u^{\prime}%
)x+1|^{2H^{\prime}-2} \nonumber\\
&&\qquad\quad{}\times \bigl[ a(H)^{2} \bigl( |u-v||u^{\prime}-v^{\prime
}||(v-v^{\prime})x+1| \bigr) ^{2H^{\prime}-2}\nonumber\\[-8pt]\\[-8pt]
&&\qquad\quad\hspace*{13.3pt}{} - 2a(H) \bigl( |u-v||(v-u^{\prime})x+1| \bigr) ^{2H^{\prime}%
-2}\nonumber\\
&&\hspace*{141.8pt}{} + |(u-u^{\prime})x+1|^{2H^{\prime}-2} \bigr]
.\nonumber
\end{eqnarray}
\end{proposition}

Before proving this proposition, let us take note of its consequence.
\begin{theorem}
\label{adjust} Let $(Z(t),t\in[0,1])$ be a Rosenblatt process with
self-similarity parameter $H\in(1/2,2/3)$ and let previous notation for
constants prevail. Then, the following convergence occurs in distribution:
\[
\lim_{_{N\rightarrow\infty}}\frac{\sqrt{N}}{\sqrt{e_{1,H}+f_{1,H}}} \biggl[
V_{N}(2,a)-\frac{\sqrt{c_{3,H}}}{N^{1-H}}Z(1) \biggr] =\mathcal{N}(0,1).
\]
\end{theorem}
\begin{pf}
By the considerations preceding the statement of Proposition \ref{AjustPro},
and (\ref{VNwantnormal}) in particular, we have that
\[
\sqrt{N} \biggl[ V_{N}(2,a)-\frac{\sqrt{c_{3,H}}}{N^{1-H}}Z(1) \biggr]
=\sqrt{N}T_{4}+\sqrt{N}U_{2}.
\]
Theorem \ref{geesbend} proves that $\sqrt{N}T_{4}$ converges in distribution
to a centered normal with variance $e_{1,H}$. Proposition \ref{AjustPro}
proves that $\sqrt{N}U_{2}$ converges in distribution to a centered normal
with variance $f_{1,H}$. Since these two sequences of random variables
live in
two distinct chaoses (fourth and second respectively), Theorem 1 in
\cite{PT}
implies that the sum of these two sequences converges in distribution
to a
centered normal with variance $e_{1,H}+f_{1,H}$. The theorem is proved.
\end{pf}

To prove Proposition \ref{AjustPro}, we must first perform the calculation
which yields the constant $f_{1,H}$ therein. This result is postponed
to the
\hyperref[Append]{Appendix}, as Lemma \ref{lemmaf1}; it shows that
$\mathbf{E}[(\sqrt{N}%
U_{2})^{2}]$ converges to $f_{1,H}$. Another (very) technical result needed
for the proof of Proposition \ref{AjustPro}, which is used to guarantee that
$\sqrt{N}U_{2}$ has a normal limiting distribution, is also included in the
\hyperref[Append]{Appendix} as Lemma \ref{LemmaAjustProof}. An
explanation of why the conclusions
of Proposition \ref{AjustPro} and Theorem \ref{adjust} cannot hold when
$H\geq2/3$ is also given in the \hyperref[Append]{Appendix}, after the
proof of Lemma \ref{LemmaAjustProof}. We now prove the proposition.
\begin{pf*}{Proof of Proposition \protect\ref{AjustPro}}
Since $U_{2}$ is a member of the second
chaos, we introduce notation for its kernel. We write%
\[
\frac{\sqrt{N}}{\sqrt{f_{1,H}}}U_{2}=I_{2} ( g_{N} ),
\]
where $g_{N}$ is the following symmetric function in $L^{2} (
[0,1]^{2} ) $:%
\[
g_{N} ( y_{1},y_{2} ) :=\frac{N^{H-1/2}}{\sqrt{f_{1,H}}} \biggl(
\frac{N^{1-H}}{4d(H)}f_{N}(y_{1},y_{2})-L_{1}(y_{1},y_{2}) \biggr) .
\]
Lemma \ref{lemmaf1} proves that $\mathbf{E} [ (
I_{2} ( g_{N} ) ) ^{2} ] = \Vert
g_{N} \Vert_{L^{2} ( [0,1]^{2} ) }^{2}$ converges
to $1$ as $N\rightarrow\infty$. By the result in \cite{NP} for
second-chaos sequences (see Theorem 1, point (ii) in \cite{NP}, which
is included as part (iii) of Theorem \ref{NOTcrit} herein), we have
that $I_{2} ( g_{N} ) $ will converge to a standard normal
if (and only
if)%
\[
\lim_{N\rightarrow\infty} \Vert g_{N}\otimes_{1}g_{N} \Vert
_{L^{2} ( [0,1]^{2} ) }^{2}=0,
\]
which would complete the proof of the proposition. This fact does hold if
$H<2/3$. We have included this technical and delicate calculation as Lemma
\ref{LemmaAjustProof} in the \hyperref[Append]{Appendix}. Following the
proof of this lemma is a
discussion of why the above limit cannot be $0$ when $H\geq2/3$.
\end{pf*}

\section{The estimators for the self-similarity parameter}\label{Stat}

In this section, we construct estimators for the self-similarity
exponent of a
Hermite process based on the discrete observations of the driving
process at
times $0,\frac{1}{N},\ldots,1$. It is known that the asymptotic
behavior of
the statistics $V_{N}(2,a)$ is related to the asymptotic properties of
a class
of estimators for the Hurst parameter $H$. This is mentioned in, for instance,
\cite{coeur}.

We recall the setup for how this works. Suppose that the observed
process $X$
is a Hermite process; it may be Gaussian (fractional Brownian motion) or
non-Gaussian (Rosenblatt process, or even a higher order Hermite
process). With
$a= \{ -1,+1 \} $, the $2$-variation is denoted by
%
\begin{equation}\label{SN}%
S_{N}(2,a)=\frac{1}{N}\sum_{i=1}^{N} \biggl( X\biggl(\frac{i}{N}\biggr)-X\biggl(\frac{i-1}%
{N}\biggr) \biggr) ^{2}.
\end{equation}
Recall that $\mathbf{E} [ S_{N}(2,a) ] =N^{-2H}.$ By estimating
$\mathbf{E} [ S_{N}(2,a) ] $ by $S_{N}(2,a)$, we can construct the
estimator
%
\begin{equation}\label{estimH}%
\hat{H}_{N}(2,a)=-\frac{\log S_{N}(2,a)}{2\log N},
\end{equation}
which coincides with the definition in (\ref{HN}) given at the
beginning of
this paper. To prove that this is a strongly consistent estimator for
$H$, we
begin by writing
\[
1+V_{N}(2,a)=S_{N}(2,a)N^{2H},%
\]
where $V_{N}$ is the original quantity defined in (\ref{VNdef0}), and thus
\[
\log\bigl( 1+V_{N}(2,a) \bigr) =\log S_{N}(2,a)+2H\log N=-2\bigl(\hat{H}%
_{N}(2,a)-H\bigr)\log N.
\]
Moreover, by Remark \ref{rem2}, $V_{N}(2,a)$ converges almost surely to $0$
and thus $\log( 1+V_{N}(2,a) ) =V_{N}(2,a)(1+o(1)),$ where
$o ( 1 ) $ converges to $0$ almost surely as $N\rightarrow\infty$.
Hence, we obtain
%
\begin{equation}\label{HV}%
V_{N}(2,a)=2\bigl(H-\hat{H}_{N}(2,a)\bigr) ( \log N ) \bigl(1+o(1)\bigr).
\end{equation}
Relation (\ref{HV}) means that the $V_{N}$'s behavior immediately gives the
behavior of $\hat{H}_{N}-H$.

Specifically, we can now state our convergence results. In the
Rosenblatt data
case, the renormalized error $\hat{H}_{N}-H$ does not converge to the normal
law. But, from Theorem \ref{adjust}, we can obtain an adjusted version
of this
error that converges to the normal distribution.
\begin{theorem}
\label{thmstat}Suppose that $H>\frac{1}{2}$ and that the observed
process $Z$ is a
Rosenblatt process with self-similarity parameter $H$. Then, strong consistency
holds for $\hat{H}_{N}$, that is, almost surely,
%
\begin{equation}\label{H1}%
\lim_{N\rightarrow\infty}\hat{H}_{N}(2,a)=H.
\end{equation}
In addition, we have the following convergence in $L^{2} ( \Omega)
$:%
%
\begin{equation}\label{H2}%
\lim_{N\rightarrow\infty}\frac{N^{1-H}}{2d(H)}\log( N ) \bigl(\hat
{H}_{N}(2,a)-H\bigr)=Z(1),
\end{equation}
where $Z ( 1 ) $ is the observed process at time $1$.

Moreover, if $H<2/3$, then, in distribution as $N\rightarrow\infty$, with
$c_{3,H}$, $e_{1,H}$ and $f_{1,H}$ in (\ref{c3H}), (\ref{e1H}) and
(\ref{f1H}),%
\[
\frac{\sqrt{N}}{\sqrt{e_{1,H}+f_{1,H}}} \biggl[ -2\log( N )
\bigl(\hat{H}_{N}(2,a)-H\bigr)-\frac{\sqrt{c_{3,H}}}{N^{1-H}}Z(1) \biggr]
\rightarrow\mathcal{{N}}(0,1).
\]
\end{theorem}
\begin{pf}
This follows from Theorems \ref{adjust} and \ref{T2toRose} and relation
(\ref{HV}).
\end{pf}

\begin{appendix}\label{Append}
\section*{Appendix}

\begin{lemma}
\label{lemmac1}The series $\sum_{k=1}^{\infty} ( 2k^{2H}- (
k-1 ) ^{2H}- ( k+1 ) ^{2H} ) ^{2}$ is finite if and
only if $H\in( 1/2,3/4 ) $.
\end{lemma}
\begin{pf}
Since $2k^{2H}- ( k-1 ) ^{2H}- ( k+1 ) ^{2H}%
=k^{2H}f ( \frac{1}{k} ) $, with $f(x):=2-(1-x)^{2H}-(1+x)^{2H}%
$ being asymptotically equivalent to $2H(2H-1)x^{2}$ for small~$x$, the
general term of the series is equivalent to $(2H)^{2} ( 2H-1 )
^{2}k^{4H-4}$.
\end{pf}
\begin{lemma}
\label{lemmac2}When $H\in( 3/4,1 ) $, $N^{2}\sum_{i,j=1,\ldots
,N; \vert i-j \vert\geq2} ( 2 \vert\frac{i-j}%
{N} \vert^{2H}- \break\vert\frac{i-j-1}{N} \vert^{2H}- \vert
\frac{i-j+1}{N} \vert^{2H} ) ^{2}$ converges to $H^{2} (
2H-1 ) / ( H-3/4 ) $ as $N\rightarrow\infty$.
\end{lemma}
\begin{pf} This is left to the reader. The proof can be found in
the extended version of this paper, available at
\url{http://arxiv.org/abs/0709.3896v2}.
\end{pf}
\begin{lemma}
\label{lemmaT2}For all $H>1/2$, with $I_{i}= ( \frac{i-1}{N},\frac{i}%
{N} ] $, $i=1,\ldots,N$,%
%
\setcounter{equation}{0}
\begin{eqnarray}\label{gg1}\quad
&& \lim_{N\rightarrow\infty}N^{2H}\sum_{i,j=1}^{N}\int_{I_{i}}\int
_{I_{i}}%
\int_{I_{j}}\int_{I_{j}}|u-v|^{2H^{\prime}-2}\nonumber\\
&&\hspace*{127.1pt}{}\times|u^{\prime}-v^{\prime
}|^{2H^{\prime}-2}|u-u^{\prime}|^{2H^{\prime}-2}\nonumber\\[-8pt]\\[-8pt]
&&\hspace*{127.1pt}{}\times|v-v^{\prime
}|^{2H^{\prime}%
-2}\,du^{\prime}\,dv^{\prime}\,dv\,du\nonumber\\
&&\qquad =2a(H)^{-2} \biggl( \frac{1}{2H-1}-\frac{1}{2H} \biggr). \nonumber
\end{eqnarray}
\end{lemma}
\begin{pf} We again refer to the extended version of the paper, online at
\href{http://arxiv.org/abs/0709.3896v2}{http://\break arxiv.org/abs/0709.3896v2}, for this proof.
\end{pf}
\begin{lemma}
\label{lemmaf1}With $f_{1,H}$ given in (\ref{f1H}) and $U_{2}$ in
(\ref{VNwantnormal}), we have
\[
\lim_{N\rightarrow\infty}\mathbf{E} \bigl[
\bigl( \sqrt{N}U_{2} \bigr) ^{2} \bigr] =f_{1,H}.
\]
\end{lemma}
\begin{pf}
We have seen that $\sqrt{c_{3,H}}=4d(H)$. We have also defined
\[
\sqrt{N}U_{2}=N^{H-1/2}\sqrt{c_{3,H}} \biggl[ \frac{N^{1-H}}{\sqrt{c_{3,H}}%
}T_{2}-Z(1) \biggr] .
\]
Let us simply compute the $L^{2}$-norm of the term in brackets. Since this
expression is a member of the second chaos and, more specifically, since
$T_{2}=I_{2} ( f_{N} ) $ and $Z ( 1 ) =I_{2} (
L_{1} ), $ where $f_{N}$ [given in (\ref{fNtrue})] and $L_{1}$ [given in
(\ref{lt})] are symmetric functions in $L^{2}([0,1]^{2})$, it holds that
\begin{eqnarray*}
&&
\mathbf{E} \biggl[ \biggl( \frac{N^{1-H}}{\sqrt{c_{3,H}}}T_{2}-Z(1) \biggr)
^{2} \biggr]\\
&&\qquad = \biggl\Vert\frac{N^{1-H}}{4d(H)}f_{N}-L_{1} \biggr\Vert
_{L^{2}([0,1]^{2})}^{2}\\
&&\qquad = \frac{N^{2-2H}}{4d(H)^{2}}\Vert f_{N}\Vert_{L^{2}([0,1]^{2})}
\\
&&\qquad\quad{}-2\frac{N^{1-H}}{4d(H)}\langle f_{N},L_{1}\rangle
_{L^{2}([0,1]^{2})} +\Vert
L_{1}\Vert_{L^{2}([0,1]^{2})}^{2}.\vadjust{\goodbreak}
\end{eqnarray*}
The first term has already been computed. It gives%
\begin{eqnarray*}
&&\frac{N^{2-2H}}{4d(H)^{2}}\Vert f_{N}\Vert_{L^{2}([0,1]^{2})}
\\
&&\qquad= N^{-2H}a^{4}(H)d^{2}(H)
\\
&&\qquad\quad{}\times
\sum_{i,j=1}^{N}\int_{[0,1]^{4}}du\,dv\,du^{\prime
}\,dv^{\prime}\\
&&\hspace*{64.7pt}{}\times ( |u-v||u^{\prime}-v^{\prime}||u-u^{\prime
}+i-j||v-v^{\prime
}+i-j| ) ^{2H^{\prime}-2}.
\end{eqnarray*}
By using the expression for the kernel $L_{1}$ and Fubini's theorem,
the scalar
product of $f_{N}$ and $L_{1}$ gives
\begin{eqnarray*}
&& \frac{N^{1-H}}{4d(H)}\langle f_{N},L_{1}\rangle_{L^{2}([0,1]^{2})}\\
&&\qquad =\int_{0}^{1}\int_{0}^{1}dy_{1}\,dy_{2}\,\frac{N^{1-H}}{4d(H)}f_{N}(y_{1}%
,y_{2})L_{1}(y_{1},y_{2})\\
&&\qquad =N^{H}a(H)^{3}d(H)^{2}\sum_{i=1}^{N}\int_{I_{i}}\int_{I_{i}}du\,dv\int
_{0}^{1}du^{\prime} ( |u-v||u-u^{\prime}||v-u^{\prime}| )
^{2H^{\prime}-2}\\
&&\qquad =N^{H}a(H)^{3}d(H)^{2}\sum_{i,j=1}^{N}\int_{I_{i}}\int_{I_{i}}%
du\,dv\int_{I_{j}}du^{\prime} ( |u-v||u-u^{\prime}||v-u^{\prime}| )
^{2H^{\prime}-2}\\
&&\qquad =N^{-2H}a(H)^{3}d(H)^{2}
\\
&&\qquad\quad{}\times\sum_{i,j=1}^{N}\int_{[0,1]^{3}} (
|u-v||u-u^{\prime}+i-j||v-u^{\prime}+i-j| ) ^{2H^{\prime}-2}%
\,du\,dv\,du^{\prime}.
\end{eqnarray*}
Finally, the last term $\Vert L_{1}\Vert_{L^{2}([0,1]^{2})}^{2}$ can be
written in the following way:
\begin{eqnarray*}
\Vert L_{1}\Vert_{L^{2}([0,1]^{2})}^{2} &=& d(H)^{2}a(H)^{2}\int
_{[0,1]^{2}%
}|u-u^{\prime}|^{2(2H^{\prime}-2)}\,du\,du^{\prime}\\
&=& d(H)^{2}a(H)^{2}\sum_{i,j=1}^{N}\int_{I_{i}}\int_{I_{j}}|u-u^{\prime
}|^{2(2H^{\prime}-2)}\,du\,du^{\prime}\\
&=& d(H)^{2}a(H)^{2}N^{-2H}\sum_{i,j=1}^{N}\int_{[0,1]^{2}}|u-u^{\prime
}+i-j|^{2(2H^{\prime}-2)}\,du\,du^{\prime}.
\end{eqnarray*}

One can check that, when bringing these three contributions together, the
``diagonal'' terms corresponding to $i=j$
vanish. Thus, we get%
\[
\mathbf{E} \bigl[ \bigl( \sqrt{N}U_{2} \bigr) ^{2} \bigr]
=32d(H)^{4}a(H)^{2}\frac{1}{N}\sum_{k=1}^{N-1}(N-k-1)k^{2H-2}F
\biggl(\frac{1}{k} \biggr),
\]
where $F$ is the function we introduced in (\ref{eff}).

This function $F$ is of class $C^{1}$ on the interval $[0,1]$. It can
be seen
that
\begin{eqnarray*}
F(0) &=& \int_{[0,1]^{4}}du\,dv\,du^{\prime}\,dv^{\prime} \\
&&\hspace*{0pt}{}\times\bigl( a(H)^{2} (
|u-v||u^{\prime}-v^{\prime}| ) ^{2H^{\prime}-2}-2a(H)|u-v|+1 \bigr) \\
&=& a(H)^{2} \biggl( \int_{[0,1]^{2}}|u-v|^{2H^{\prime}-2} \biggr)
^{2}-2a(H)\int_{[0,1]^{2}}|u-v|^{2H^{\prime}-2}\,du\,dv+1\\
&=& 0.
\end{eqnarray*}
Similarly, one can also calculate the derivative $F^{\prime}$ and check that
$F^{\prime} ( 0 ) =0$. Therefore, $F(x)=o(x)$ as $x\rightarrow0$. To
investigate the sequence $a_{N}:=N^{-1}\sum_{k=1}^{N-1}(N-k-1)k^{2H-2}F (
\frac{1}{k} ) $, we split it into two pieces:%
\begin{eqnarray*}
a_{N} &=& N^{-1}\sum_{k=1}^{N-1}(N-k-1)k^{2H-2}F \biggl( \frac{1}{k} \biggr) \\
&=& \sum_{k=1}^{N-1}k^{2H-2}F \biggl( \frac{1}{k} \biggr) +N^{-1}\sum
_{k=1}^{N-1}(k+1)k^{2H-2}F \biggl( \frac{1}{k} \biggr) \\
&=:& b_{N}+c_{N}.%
\end{eqnarray*}
Since $b_{N}$ is the partial sum of a sequence of positive terms, one only
needs to check that the series is finite. The relation $F ( 1/k )
\ll1/k$ yields that it is finite if and only if $2H-3<-1$, which is
true. For the term
$c_{N}$, one notes that we may replace the factor $k+1$ by $k$ since,
by the
calculation undertaken for $b_{N}$, $N^{-1}\sum_{k=1}^{N-1}k^{2H-2}F (
\frac
{1}{k} ) $ converges to $0$. Hence, asymptotically, we have%
\[
c_{N}\simeq N^{-1}\sum_{k=1}^{N-1}k^{2H-3}F \biggl( \frac{1}{k} \biggr) \leq
N^{-1} \Vert F \Vert_{\infty}\sum_{k=1}^{\infty}k^{2H-3},
\]
which thus converges to $0$. We have proven that $\lim a_{N}=\lim b_{N}%
=\sum_{k=1}^{\infty}k^{2H-2}\times\break F ( \frac{1}{k} ) $, which completes the
proof of the lemma.
\end{pf}
\begin{lemma}
\label{LemmaAjustProof}
Defining
\[
g_{N} ( y_{1},y_{2} ) :=\frac{N^{H-1/2}}{\sqrt{f_{1,H}}}
\biggl(\frac{N^{1-H}}{4d(H)}f_{N}(y_{1},y_{2})-L_{1}(y_{1},y_{2}) \biggr),\vadjust{\goodbreak}
\]
we have $\lim_{N\rightarrow\infty} \Vert g_{N}\otimes_{1}g_{N} \Vert
_{L^{2} ( [0,1]^{2} ) }^{2}=0$ provided $H<2/3$.
\end{lemma}
\begin{pf}
We omit the leading constant $f_{1,H}^{-1/2}$, which is irrelevant.
Using the
expression (\ref{fNtrue}) for $f_{N}$, we have%
\begin{eqnarray*}
g_{N}(y_{1},y_{2})&=&N^{2H-1/2}d(H)a(H)\\
&&{}\times\sum_{i=1}^{N}\int_{I_{i}}\int
_{I_{i}%
}\partial_{1}K^{H^{\prime}}(u,y_{1})\,\partial_{1}K^{H^{\prime}}(v,y_{2}%
)|u-v|^{2H^{\prime}-2}\,dv\,du\\
&&{} -L_{1}(y_{1},y_{2}).
\end{eqnarray*}
Here, and below, we will be omitting indicator functions of the
type\break
$1_{[0,({i+1})/{N}]}(y_{1})$ because, as stated earlier, these are
implicitly contained in the support of $\partial_{1}K^{H^{\prime}}$. By
decomposing the expression for $L_{1}$ from (\ref{lt}) over the same blocks
$I_{i}\times I_{i}$ as for $f_{N}$, we can now express the contraction
$g_{N}\otimes_{1}g_{N}$ as follows:
\[
(g_{N}\otimes_{1}g_{N})(y_{1},y_{2})=N^{2H-1} ( A_{N}-2B_{N}%
+C_{N} ),
\]
where we have introduced three new quantities,%
\begin{eqnarray*}
A_{N} &:=& N^{2H}d(H)^{2}a(H)^{3}\\
&&{}\times\sum_{i,j=1}^{N}\int_{I_{i}}\int_{I_{i}%
}dv\,du\int_{I_{j}}\int_{I_{j}}dv^{\prime}\,du^{\prime}\\
&&\hspace*{0pt}{} \times[ |u-v|\cdot|u^{\prime}-v^{\prime}|\cdot|v-v^{\prime}| ]
^{2H^{\prime}-2}\\
&&\hspace*{0pt}{} \times\partial_{1}K^{H^{\prime}}(u,y_{1})\partial
_{1}K^{H^{\prime}%
}(u^{\prime},y_{2}),
\\
B_{N} &:=& N^{H}a(H)^{2}d(H)^{2}\sum_{i=1}^{N}\int_{I_{i}}\int_{I_{i}}%
dv\,du\int_{0}^{1}du^{\prime} [ |u-v|\cdot|u^{\prime}-v| ]
^{2H^{\prime}-2}\\
&&\hspace*{0pt}{} \times\partial_{1}K^{H^{\prime}}(u,y_{1})\,\partial
_{1}K^{H^{\prime}%
}(u^{\prime},y_{2})\\
&=& N^{H}a(H)^{2}d(H)^{2}\sum_{i,j=1}^{N}\int_{I_{i}}\int_{I_{i}}%
dv\,du\int_{I_{j}}du^{\prime} [ |u-v|\cdot|u^{\prime}-v| ]
^{2H^{\prime}-2}\\
&&\hspace*{0pt}{}\times \partial_{1}K^{H^{\prime}}(u,y_{1})\,\partial
_{1}K^{H^{\prime}%
}(u^{\prime},y_{2})
\end{eqnarray*}
and
\begin{eqnarray*}
C_{N} &=& d(H)^{2}a(H)\int_{0}^{1}\int_{0}^{1}dv\,du\,\partial
_{1}K^{H^{\prime}%
}(u,y_{1})\,\partial_{1}K^{H^{\prime}}(v,y_{2})|u-v|^{2H^{\prime}-2}\\
&=& d(H)^{2}a(H)\sum_{i,j=1}^{N}\int_{I_{i}}\int_{I_{j}}dv\,du\,\partial
_{1}K^{H^{\prime}}(u,y_{1})\,\partial_{1}K^{H^{\prime}}(v,y_{2}%
)|u-v|^{2H^{\prime}-2}.\vadjust{\goodbreak}
\end{eqnarray*}
The squared norm of the contraction can then be written as%
\begin{eqnarray*}
&&\Vert g_{N}\otimes_{1}g_{N}\Vert_{L^{2}([0,1]^{2})}^{2} \\
&&\qquad =N^{4H-2} \bigl(
\Vert A_{N}\Vert_{L^{2}([0,1]^{2})}^{2}+4\Vert B_{N}\Vert
_{L^{2}([0,1]^{2}%
)}^{2}\\
&&\qquad\quad\hspace*{35.1pt}{} + \Vert C_{N}\Vert_{L^{2}([0,1]^{2})}^{2}
-4\langle A_{N},B_{N}\rangle_{L^{2}([0,1]^{2})}\\
&&\qquad\quad\hspace*{35.1pt}{} + 2\langle
A_{N},C_{N}\rangle_{L^{2}([0,1]^{2})}-4\langle B_{N},C_{N}\rangle
_{L^{2}([0,1]^{2})} \bigr) .
\end{eqnarray*}
Using the definitions of $A_{N}$, $B_{N}$ and $C_{N}$, we may express
all six
terms above explicitly. All of the computations are based on the key relation
(\ref{souvent}).

We obtain
\begin{eqnarray*}
&& \Vert A_{N}\Vert_{L^{2}([0,1]^{2})}^{2}\\[0.3pt]
&&\qquad =N^{4H}a(H)^{6}d(H)^{4}a(H)^{2}\\[0.3pt]
&&\qquad\quad{}\times\sum_{i,j,k,l=1}^{N}\int_{I_{i}}\int
_{I_{i}%
}dv\,du\int_{I_{j}}\int_{I_{j}}dv^{\prime}\,du^{\prime}\\[0.3pt]
&&\qquad\quad{}\times\int_{I_{k}}\int
_{I_{k}%
}d\bar{u}\,d\bar{v}\int_{I_{l}}\int_{I_{l}}d\bar{u}^{\prime}\,d\bar{v}^{\prime}\\[0.3pt]
&&\qquad\quad{}\times [ |u-v|\cdot|u^{\prime}-v^{\prime}|\cdot|v-v^{\prime}|\cdot|\bar
{u}-\bar{v}|\cdot|\bar{u}^{\prime}-\bar{v}^{\prime}|\\
&&\hspace*{118.5pt}{}\times|\bar{v}%
-\bar{v}^{\prime}|\cdot|u-\bar{u}|\cdot|u^{\prime}-\bar{u}^{\prime}| ]
^{2H^{\prime}-2}\\[0.3pt]
&&\qquad
=N^{4H}a(H)^{8}d(H)^{4}\frac{1}{N^{8}}\frac{1}{N^{8(2H^{\prime}-2)}}\\
&&\qquad\quad{}\times
\sum_{i,j,k,l=1}^{N}\int_{[0,1]^{8}}du\,dv\,du^{\prime}\,dv^{\prime}\,d\bar
{u}\,d\bar
{v}\,d\bar{u}^{\prime}\,d\bar{v}^{\prime}\\[0.3pt]
&&\qquad\quad\hspace*{72.5pt}\hspace*{-29.1pt}{}\times \bigl\vert|u-v|\cdot|u^{\prime}-v^{\prime}||\bar{u}-\bar{v}||\bar
{u^{\prime}}-\bar{v}^{\prime}| \bigr\vert^{2H^{\prime}-2}\\[0.3pt]
&&\qquad\quad\hspace*{72.5pt}\hspace*{-29.1pt}{}\times [
|v-v^{\prime}+i-j|\cdot|\bar{v}-\bar{v}^{\prime}+k-l|\\[0.3pt]
&&\qquad\quad\hspace*{72.5pt}\hspace*{-29.1pt}\hspace*{14.5pt}{}\times
|u-\bar{u}+i-k|\cdot|u^{\prime}-\bar{u}^{\prime}+j-l| ] ^{2H^{\prime}%
-2},
\\[0.3pt]
&&\Vert B_{N}\Vert_{L^{2}([0,1]^{2})}^{2}\\[0.3pt]
&&\qquad =N^{2H}a(H)^{6}d(H)^{4}\\[0.3pt]
&&\qquad\quad{}\times
\sum_{i,j,k,l=1}^{N}\int_{I_{i}}\int_{I_{i}}dv\,du\int_{I_{j}}du^{\prime}%
\int_{I_{k}}\int_{I_{k}}d\bar{u}\,d\bar{v}\int_{I_{l}}d\bar{u}^{\prime}\\[0.3pt]
&&\qquad\quad{}\times [ |u-v|\cdot|u^{\prime}-v||\bar{u}-\bar{v}|\cdot|\bar{u}^{\prime}%
-\bar{v}|\cdot|u-\bar{u}|\cdot|u^{\prime}-\bar{u}^{\prime}| ]
^{2H^{\prime}-2}
\end{eqnarray*}
\begin{eqnarray*}
&&\qquad =N^{2H}a(H)^{6}d(H)^{4}\\
&&\qquad\quad{}\times\sum_{i,j,k,l=1}^{N}\int
_{[0,1]^{6}}du\,dv\,du^{\prime
}\,d\bar{u}\,d\bar{v}\,d\bar{u}^{\prime}\\
&&\qquad\quad\hspace*{44.2pt}{}\times [ |u-v|\cdot|u^{\prime}-v+i-j||\bar{u}-\bar{v}|\cdot|\bar{u}^{\prime}%
-\bar{v}+k-l|\\
&&\qquad\quad\hspace*{116.3pt}{}\times|u-\bar{u}+i-k|\cdot|u^{\prime}-\bar{u}^{\prime}%
+j-l| ] ^{2H^{\prime}-2}%
\end{eqnarray*}
and
\begin{eqnarray*}
&& \Vert C_{N}\Vert_{L^{2}([0,1]^{2})}^{2}\\
&&\qquad =N^{2H}a(H)^{4}d(H)^{4}\sum_{i,j,k,l=1}^{N}\int_{I_{i}}\int_{I_{j}}%
dv\,du\int_{I_{k}}\int_{I_{l}}dv^{\prime}\,du^{\prime} \\
&&\qquad\quad{}\times[ |u-v|\cdot
|u^{\prime}-v^{\prime}|\cdot|u-u^{\prime}|\cdot|v-v^{\prime}| ]
^{2H^{\prime}-2}\\
&&\qquad
=N^{2H}a(H)^{4}d(H)^{4}\frac{1}{N^{4}}\frac{1}{N^{4(2H^{\prime}-2)}}\\
&&\qquad\quad{}\times\sum_{i,j,k,l=1}^{N}\int_{[0,1]^{4}}du\,dv\,du^{\prime}\,dv^{\prime}\\
&&\qquad\quad{}\times [ |u-v+i-j|\cdot|u^{\prime}-v^{\prime}+k-l|\\
&&\qquad\quad\hspace*{14.51pt}{}\times|u-u^{\prime
}+i-k|\cdot|v-v^{\prime}+j-l| ] ^{2H^{\prime}-2}.%
\end{eqnarray*}
The inner product terms can be also treated in the same manner. First,
\begin{eqnarray*}
&&\hspace*{-8pt}\langle A_{N},B_{N}\rangle_{L^{2}([0,1]^{2})}\\
&&\hspace*{-8pt}\qquad =N^{3H}a(H)^{7}d(H)^{4}\\
&&\hspace*{-8pt}\qquad\quad{}\times\sum_{i,j,k,l=1}^{N}\int_{I_{i}}\int_{I_{i}}%
du\,dv\int_{I_{j}}\int_{I_{j}}du^{\prime}\,dv^{\prime}\int_{I_{k}}\int
_{I_{k}%
}d\bar{u}\,d\bar{v}\int_{I_{l}}d\bar{u}^{\prime}\\
&&\hspace*{-8pt}\qquad\quad{}\times [ |u-v|\cdot|u^{\prime}-v^{\prime}|\cdot|v-v^{\prime}|\cdot|\bar
{u}-\bar{v}|\cdot|\bar{u}^{\prime}-\bar{v}|\cdot|u-\bar{u}|\cdot
|u^{\prime
}-\bar{u}^{\prime}| ] ^{2H^{\prime}-2}\\
&&\hspace*{-8pt}\qquad =N^{3H}a(H)^{7}d(H)^{4}\frac{1}{N^{7}}\frac{1}{N^{7(2H^{\prime}-2)}}
\\
&&\hspace*{-8pt}\qquad\quad{}\times\sum_{i,j,k,l=1}^{N}\int_{[0,1]^{7}}du\,dv\,du^{\prime}\,dv^{\prime}\,d\bar
{u}\,d\bar
{v}\,d\bar{u}^{\prime}\\
&&\hspace*{-8pt}\qquad\quad\hspace*{44.4pt}{}\times [ |u-v|\cdot|u^{\prime}-v^{\prime}|\cdot|v-v^{\prime}+i-j|\cdot
|\bar{u}-\bar{v}|\\
&&\hspace*{-8pt}\qquad\quad\hspace*{59.7pt}{}\times|\bar{u}^{\prime}-\bar{v}+k-l|\cdot|u-\bar{u}
+i-k|\cdot|u^{\prime}-\bar{u}^{\prime}+j-l| ]
^{2H^{\prime}-2}
\end{eqnarray*}
and
\begin{eqnarray*}
&&\langle A_{N},C_{N}\rangle_{L^{2}([0,1]^{2})}\\
&&\qquad =N^{2H}a(H)^{6}d(H)^{4}\\
&&\qquad\quad{}\times\sum_{i,j,k,l=1}^{N}\int_{I_{i}}\int_{I_{i}}du\,dv\int_{I_{j}}\int_{I_{j}%
}du^{\prime}\,dv^{\prime}\int_{I_{k}}d\bar{u}\int_{I_{l}}d\bar{v}\\
&&\qquad\quad{}\times [ |u-v|\cdot|u^{\prime}-v^{\prime}|\cdot|v-v^{\prime}|\cdot|\bar
{u}-\bar{v}|\cdot|u-\bar{u}|\cdot|u^{\prime}-\bar{v} ] ^{2H^{\prime}%
-2}\\
&&\qquad
=N^{2H}a(H)^{6}d(H)^{4}\frac{1}{N^{6}}\frac{1}{N^{6(2H^{\prime}-2)}}\\
&&\qquad\quad{}\times
\sum_{i,j,k,l=1}^{N}\int_{[0,1]^{6}}du\,dv\,du^{\prime}\,dv^{\prime}\,d\bar
{u}\,d\bar
{v}\\
&&\qquad\quad\hspace*{44.5pt}{}\times [
|u-v|\cdot|u^{\prime}-v^{\prime}|\cdot|v-v^{\prime}+i-j|\\
&&\hspace*{92.6pt}{}\times
|u-\bar{u}+i-k|\cdot|\bar{u}-\bar{v}+k-l|\cdot u^{\prime}-\bar{v} ]
^{2H^{\prime}-2}.%
\end{eqnarray*}
Finally,
\begin{eqnarray*}
&&\langle B_{N},C_{N}\rangle_{L^{2}([0,1]^{2})} \\
&&\qquad =N^{H}a(H)^{3}d(H)^{4}\\
&&\qquad\quad{}\times\sum_{i,j,k,l=1}^{N}\int_{I_{i}}\int_{I_{i}}du\,dv\int_{I_{j}}du^{\prime}%
\int_{I_{k}}d\bar{u}\int_{I_{l}}d\bar{v}\\
&&\qquad\quad{}\times [ |u-v|\cdot|u^{\prime}-v|\cdot|\bar{u}-\bar{v}|\cdot|u-\bar{u}%
|\cdot|u^{\prime}-\bar{v}| ] ^{2H^{\prime}-2}\\
&&\qquad
=N^{H}a(H)^{3}d(H)^{4}\frac{1}{N^{5}}\frac{1}{N^{5(2H^{\prime}-2)}}\\
&&\qquad\quad{}\times\sum_{i,j,k,l=1}^{N}\int_{[0,1]^{5}}du\,dv\,du^{\prime}\,d\bar{u}\,d\bar{v}\\
&&\qquad\quad\hspace*{44.7pt}{}\times [
|u-v|\cdot|u^{\prime}-v+i-j|\cdot|\bar{u}-\bar{v}+k-l|\\
&&\qquad\quad\hspace*{86.1pt}{}\times
|u-\bar{u}+i-k|\cdot|u^{\prime}-\bar{v}+j-l| ] ^{2H^{\prime}-2}.
\end{eqnarray*}

We now summarize our computations. Note that the factors $d(H)^{4}$ and
$\frac{1}{N^{4}}\frac{1}{N^{4(2H^{\prime}-2)}}$ are common to all
terms. We
also note that any terms corresponding to difference of indices smaller than
$3$ can be shown to tend collectively to $0$, similarly for other
``diagonal'' terms in this study. The proof
is omitted. We thus assume that the sums over the set $D$ of indices $i,j,k,l$
in $\{1,\ldots,N\}$\vadjust{\goodbreak} such that $ \vert i-j \vert$, $ \vert
k-l \vert$, $ \vert i-k \vert$ and $ \vert j-l \vert$
are all at least $2$. Hence, we get
%
\begin{eqnarray}\label{gn-contra}\hspace*{22pt}
&&\Vert g_{N}\otimes_{1}g_{N} \Vert_{L^{2}([0,1]^{2})}%
^{2}\nonumber\\
&&\qquad= d(H)^{4}N^{4H-2}\frac{1}{N^{4}}\nonumber\\[-8pt]\\[-8pt]
&&\qquad\quad{}\times
\sum_{(i,j,k,l)\in D} \biggl( \frac{|i-j|\cdot|k-l|\cdot|i-k|\cdot|j-l|}{N^{4}} \biggr)
^{2H^{\prime}-2}\nonumber\\
&&\qquad\quad\hspace*{51.2pt}{}\times
G\biggl(\frac{1}{i-j},\frac{1}{k-l},\frac{1}{i-k},\frac{1}{j-l}\biggr),\nonumber
\end{eqnarray}
where the function $G$ is defined for $(x,y,z,w)\in[1/2,1/2]^{4}$ by
\begin{eqnarray*}
&& G(x,y,z,w)\\
&&\qquad
=a(H)^{8}\int_{[0,1]^{8}}du\,dv\,du^{\prime}\,dv^{\prime}\,d\bar{u}\,d\bar{v}\,
d\bar{u}^{\prime}\,d\bar{v}^{\prime}\\
&&\qquad\quad{}\times[ |u-v|\cdot|u^{\prime}-v^{\prime
}|\cdot|\bar{u}-\bar{v}|\cdot|\bar{u}^{\prime}-\bar{v}^{\prime}| ]
^{2H^{\prime}-2}\\
&&\qquad\quad{}\times [ |(v-v^{\prime})x+1|\cdot|(\bar{v}-\bar{v}^{\prime})y+1|\\
&&\qquad\quad\hspace*{14.7pt}{}\times
|(u-\bar{u})z+1|\cdot|(u^{\prime}-\bar{u}^{\prime})w+1| ] ^{2H^{\prime
}-2}\\
&&\qquad\quad{} +4a(H)^{6}\int_{[0,1]^{6}}du\,dv\,du^{\prime}\,d\bar{u}\,d\bar{v}\,d\bar
{u}^{\prime}%
\\
&&\qquad\quad\hspace*{9pt}{}\times [ |u-v|\cdot|\bar{u}-\bar{v}|\cdot|(u^{\prime}-v)x+1|\cdot
|(\bar{u}^{\prime}-\bar{v})y+1|\\
&&\qquad\quad\hspace*{84.4pt}{}\times|(u-u^{\prime})z+1|\cdot|(u^{\prime}
-\bar{u}^{\prime})w+1| ] ^{2H^{\prime}-2}\\
&&\qquad\quad{} +a(H)^{4}\int_{[0,1]^{4}}du\,dv\,du^{\prime}dv\,dv^{\prime} \\
&&\qquad\quad\hspace*{9pt}{}\times[ |(u-v)x+1|\cdot
|(u^{\prime}-v^{\prime})y+1|\\
&&\qquad\quad\hspace*{23.7pt}{}\times|(u-u^{\prime})z+1|\cdot|(v-v^{\prime
})w+1| ] ^{2H^{\prime}-2}\\
&&\qquad\quad{}
-4a(H)^{7}\int_{[0,1]^{7}}du\,dv\,du^{\prime}\,dv^{\prime}\,d\bar{u}\,d\bar{v}\,
d\bar{u}^{\prime}\\
&&\qquad\quad\hspace*{9pt}{}\times [ |u-v|\cdot|u^{\prime}-v^{\prime}|\cdot|\bar{u}-\bar{v}%
|\cdot|(v-v^{\prime})x+1|\\
&&\qquad\quad\hspace*{23.7pt}{}\times|(\bar{u}^{\prime}-\bar{v})y+1|\cdot
|(u-u^{\prime})z+1|\cdot|(u^{\prime}-\bar{u}^{\prime})w+1| ]
^{2H^{\prime}-2}\\
&&\qquad\quad{} +2a(H)^{6}\int_{[0,1]^{6}}du\,dv\,du^{\prime}\,dv^{\prime}\,d\bar{u}\,d\bar{v}\\
&&\qquad\quad\hspace*{9pt}{} \times[ |u-v|\cdot|u^{\prime}-v^{\prime}|\cdot|(v-v^{\prime})x+1|\cdot
|(\bar{u}-\bar{v})y+1|\\
&&\qquad\quad\hspace*{23.7pt}\hspace*{66.8pt}{}\times|(u-u^{\prime})z+1|\cdot|(u^{\prime}-\bar
{v})w+1| ] ^{2H^{\prime}-2}\\
&&\qquad\quad{} -4a(H)^{5}\int_{[0,1]^{5}}du\,dv\,du^{\prime}\,d\bar{u}\,d\bar{v}\\
&&\qquad\quad\hspace*{9pt}{}\times [ |u-v|\cdot|(v-u^{\prime})x+1|\cdot|(\bar{u}-\bar{v})y+1|\\
&&\qquad\quad\hspace*{49pt}{}\times
|(u-\bar{u})z+1|\cdot|(u^{\prime}-\bar{v})w+1| ] ^{2H^{\prime}-2}.
\end{eqnarray*}

It is elementary to check that $G$ and all its partial derivatives are bounded
on $[-1/2,1/2]^{4}$. More specifically, by using the identity
\[
a(H)^{-1}=\int_{0}^{1}\int_{0}^{1}|u-v|^{2H^{\prime}-2}\,du\,dv,
\]
we obtain
\begin{eqnarray*}
G(0,0,0,0) &=&
a(H)^{4}+4a(H)^{4}+a(H)^{4}-4a(H)^{4}+2a(H)^{4}-4a(H)^{4}\nonumber\\[-8pt]\\[-8pt]
&=& 0.\nonumber
\end{eqnarray*}
The boundedness of $G$'s partial derivatives implies, by the mean value
theorem, that there exists a constant $K$ such that, for all $ (
i,j,k,l ) \in D$,%
\begin{eqnarray*}
&&\biggl\vert G\biggl(\frac{1}{i-j},\frac{1}{k-l},\frac{1}{i-k},\frac{1}{j-l}%
\biggr) \biggr\vert\\
&&\qquad\leq\frac{K}{ \vert i-j \vert}+\frac{K}{ \vert
k-l \vert}+\frac{K}{ \vert i-k \vert}+\frac{K}{ \vert
j-l \vert}.
\end{eqnarray*}
Hence, from (\ref{gn-contra}), because of the symmetry of the sum with respect
to the indices, it is sufficient to show that the following converges
to $0$:%
%
\begin{equation}\label{ess}\hspace*{28pt}
S:=N^{4H-2}\frac{1}{N^{4}}\sum_{(i,j,k,l)\in D} \biggl( \frac{|i-j|\cdot
|k-l|\cdot|i-k|\cdot|j-l|}{N^{4}} \biggr) ^{H-1}\frac{1}{ \vert
i-j \vert}.
\end{equation}
We will express this quantity by singling out the term $i^{\prime
}:=i-j$ and
summing over it last:%
\begin{eqnarray*}
S
&=& 2N^{4H-1}\sum_{i^{\prime}=3}^{N-1}\frac{1}{N^{3}}\mathop{\mathop{\sum}
_{(i^{\prime}+j,j,k,l)\in D}}_{1\leq j\leq N-i^{\prime}} \biggl(
\frac{|k-l|\cdot|i^{\prime}+j-k|\cdot|j-l|}{N^{3}} \biggr) ^{H-1} \biggl(
\frac{i^{\prime}}{N} \biggr) ^{H-1}\frac{1}{i^{\prime}}\\
&=& 2N^{3H-2}\sum_{i^{\prime}=3}^{N-1} ( i^{\prime} ) ^{H-2}%
\frac{1}{N^{3}}\mathop{\mathop{\sum}_{(i^{\prime}+j,j,k,l)\in
D}}_{1\leq j\leq
N-i^{\prime}} \biggl( \frac{|k-l|\cdot|i^{\prime}+j-k|\cdot|j-l|}{N^{3}%
} \biggr) ^{H-1}.
\end{eqnarray*}
For fixed $i^{\prime}$, we can compare the sum over $j,k,l$ to a Riemann
integral since the power $H-1>-1$. This cannot be done, however, for $ (
i^{\prime} ) ^{H-2}$; rather, one must use the fact that this is the
term of a summable series. We get that, asymptotically for large $N$,%
\[
S\simeq2N^{3H-2}\sum_{i^{\prime}=3}^{N-1} ( i^{\prime} )
^{H-2}g ( i^{\prime}/N ),
\]
where the function $g$ is defined on $[0,1]$ by%
%
\begin{equation}\label{nuthergee}%
g ( x ) :=\int_{0}^{1-x}\int_{0}^{1}\int_{0}^{1}dy\,dz\,dw\, \vert
( z-w ) ( x+y-z ) ( y-w ) \vert
^{H-1}.
\end{equation}
It is easy to check that $g$ is a bounded function on $[0,1]$; thus, we have
proven that for some constant $K>0$,%
\[
S\leq KN^{3H-2}\sum_{i^{\prime}=3}^{\infty} ( i^{\prime} ) ^{H-2},%
\]
which converges to $0$ provided $H<2/3$. This completes the proof of
the lemma.
\end{pf}

We conclude this appendix with a discussion of why the threshold
$H<2/3$ cannot
be improved upon, and the consequences of this. We can perform a finer
analysis of the function $G$ in the proof above. The first and second
derivatives of $G$ at $\bar{0}= ( 0,0,0,0 ) $ can be calculated by
hand. The calculation is identical for $\partial G/\partial x ( \bar
{0} ) $ and all other first derivatives, yielding [via the
expression used above for $a ( H ) $],%
\begin{eqnarray*}
&& \frac{1}{H-1}\frac{\partial G}{\partial x} ( \bar{0} ) \\
&&\qquad =a(H)^{6}\int_{[0,1]^{4}}du\,dv\,du^{\prime}\,dv^{\prime} ( v-v^{\prime
} ) [ |u-v|\cdot|u^{\prime}-v^{\prime}| ] ^{H-1}\\
&&\qquad\quad{} +4a(H)^{5}\int_{[0,1]^{3}}du\,dv\,du^{\prime} ( v-u^{\prime} )
|u-v|^{H-1}\\
&&\qquad\quad{} +a(H)^{4}\int_{[0,1]^{2}}du\,dv ( u-v )\\
&&\qquad\quad{} -4a(H)^{6}\int_{[0,1]^{4}%
}du\,dv\,du^{\prime}\,dv^{\prime} ( v-v^{\prime} ) [ |u-v|\cdot
|u^{\prime}-v^{\prime}| ] ^{H-1}\\
&&\qquad\quad{} +2a(H)^{6}\int_{[0,1]^{4}}du\,dv\,du^{\prime}\,dv^{\prime} ( v-v^{\prime
} ) [ |u-v|\cdot|u^{\prime}-v^{\prime}| ] ^{H-1}\\
&&\qquad\quad{}-4a(H)^{5}\int_{[0,1]^{3}}du\,dv\,du^{\prime} ( v-u^{\prime} )
\vert u-v \vert^{H-1}.
\end{eqnarray*}
We note that the two lines with $4a ( H ) ^{5}$ cancel each other
out. For each of the other four lines, we see that the factor $ (
v-v^{\prime} ) $ is an odd term and the other factor is symmetric with
respect to $v$ and $v^{\prime}$. Therefore, each of the other four
factors is
zero individually. This proves that the gradient of $G$ at $0$ is null.
Let us
find expressions for the second derivatives. Similarly to the
above calculation, we can
write%
\begin{eqnarray*}
&& \frac{1}{ ( 1-H ) ( 2-H ) }\frac{\partial^{2}%
G}{\partial x^{2}} ( \bar{0} ) \\
&&\qquad =a(H)^{6}\int_{[0,1]^{4}}du\,dv\,du^{\prime}\,dv^{\prime} ( v-v^{\prime
} ) ^{2} [ |u-v|\cdot|u^{\prime}-v^{\prime}| ] ^{H-1}\\
&&\qquad\quad{}
+4a(H)^{5}\int_{[0,1]^{3}}du\,dv\,du^{\prime} ( v-u^{\prime} )
^{2}|u-v|^{H-1}\\
&&\qquad\quad{} +a(H)^{4}\int_{[0,1]^{2}}du\,dv ( u-v ) ^{2}\\
&&\qquad\quad{} -4a(H)^{6}%
\int_{[0,1]^{4}}du\,dv\,du^{\prime}\,dv^{\prime} ( v-v^{\prime} )
^{2} [ |u-v|\cdot|u^{\prime}-v^{\prime}| ] ^{H-1}\\
&&\qquad\quad{} +2a(H)^{6}\int_{[0,1]^{4}}du\,dv\,du^{\prime}\,dv^{\prime} ( v-v^{\prime
} ) ^{2} [ |u-v|\cdot|u^{\prime}-v^{\prime}| ] ^{H-1}\\
&&\qquad\quad{}
-4a(H)^{5}\int_{[0,1]^{3}}du\,dv\,du^{\prime} ( v-u^{\prime} )
^{2} \vert u-v \vert^{H-1}.
\end{eqnarray*}
Again, the terms with $a ( H ) ^{5}$ cancel each other out. The
three terms with $a ( H ) ^{6}$ add to a nonzero value and we
thus get
\begin{eqnarray*}
&&\frac{1}{ ( 1-H ) ( 2-H ) }\frac{\partial^{2}
G}{\partial x^{2}} ( \bar{0} ) \\
&&\qquad =-a(H)^{6}\int_{[0,1]^{4}%
}du\,dv\,du^{\prime}\,dv^{\prime} ( v-v^{\prime} ) ^{2} [
|u-v|\cdot|u^{\prime}-v^{\prime}| ] ^{H-1}\\
&&\qquad\quad{} +a(H)^{4}\int_{[0,1]^{4}}du\,dv ( u-v ) ^{2}.
\end{eqnarray*}
While the evaluation of this integral is nontrivial, we can show that for
all $H>1/2$, it is a strictly positive constant $\gamma( H ) $.
Similar computations can be attempted for the mixed derivatives, which
are all
equal to some common value $\eta( H ) $ at $\bar{0}$ because of
$G$'s symmetry, and we will see that the sign of $\eta( H ) $ is
irrelevant. We can now write, using Taylor's formula,%
\begin{eqnarray*}
G ( x,y,z,w ) & = & \gamma( H ) ( x^{2}%
+y^{2}+z^{2}+w^{2} ) \\
&&{} +\eta( H ) ( xy+xz+xw+yz+yw+zw )\\
&&{} +o (
x^{2}+y^{2}+z^{2}+w^{2} ) .
\end{eqnarray*}
By taking $x^{2}+y^{2}+z^{2}+w^{2}$ sufficiently small [this
corresponds to
restricting $ \vert i-j \vert$ and other differences to being larger
than some value $m=m ( H ) $, whose corresponding ``diagonal''
terms not satisfying this restriction are dealt with as
usual], we get, for some constant $\theta( H ) >0$,%
\[
G ( x,y,z,w ) \geq\theta( H ) ( x^{2}+y^{2}%
+z^{2}+w^{2} ) +\eta( H ) ( xy+xz+xw+yz+yw+zw )
.
\]

Let us first look at the terms in (\ref{gn-contra}) corresponding to
$x^{2}+y^{2}+z^{2}+w^{2}$. These are collectively bounded below by the same
sum restricted to $i=j+m$, which equals%
\[
d(H)^{4}N^{4H-2}\frac{1}{N^{4}}\sum_{(j+m,j,k,l)\in D} \biggl( \frac
{|i-j|\cdot|k-l|\cdot|i-k|\cdot|j-l|}{N^{4}} \biggr) ^{2H^{\prime}-2}%
\frac{\theta( H ) }{ ( i-j ) ^{2}}.
\]
The fact that the final factor contains $ ( i-j ) ^{-2}$ instead of
$ ( i-j )^{-1}$, which we had, for instance, in (\ref{ess}) in the
proof of the lemma, does not help us. In particular, calculations
identical to those following (\ref{ess}) show that the above is larger than
\[
2N^{3H-2}g ( m/N ),
\]
which does not go to $0$ if $H\geq2/3$ since $g ( 0 ) $ calculated
from (\ref{nuthergee}) is positive.

For the terms in (\ref{gn-contra}) corresponding to $xy+xz+xw+yz+yw+zw$,
considering, for instance, the term $xy$, similar computations to those above
lead to the corresponding term in $S$ being equal to
\begin{eqnarray*}
&& 2N^{2H-2}\sum_{i^{\prime}=m}^{N-1}\sum_{k^{\prime}=m}^{N-1} (
i^{\prime}k^{\prime} ) ^{H-2}\frac{1}{N^{2}}\\
&&\quad\mathop{\mathop{\sum
}_{(i^{\prime
}+j,j,k^{\prime}+l,l)\in D}}_{1\leq j\leq N-i^{\prime};1\leq l\leq
N-k^{\prime}%
} \biggl( \frac{|i^{\prime}+j-k^{\prime}-l|\cdot|j-l|}{N^{3}} \biggr) ^{H-1}\\
&&\qquad \simeq2N^{2H-2}\sum_{i^{\prime}=m}^{N-1}\sum_{k^{\prime}=m}^{N-1} (
i^{\prime}k^{\prime} ) ^{H-2}\int_{0}^{1-i^{\prime}/N}\int
_{0}^{1-k^{\prime}/N}dy\,dw \\
&&\hspace*{109.6pt}{}\times\biggl\vert ( z-w ) \biggl( \frac
{i^{\prime}}{N}+y-\frac{k^{\prime}}{N}-w \biggr) ( y-w )
\biggr\vert^{H-1},
\end{eqnarray*}
which evidently tends to $0$ provided $H<1$.

We conclude that if $H\geq2/3$, then $ \Vert g_{N}\otimes_{1}g_{N} \Vert
_{L^{2}([0,1]^{2})}^{2}$ does not tend to $0$ and, by the
Nualart--Ortiz-Latorre criterion [Theorem \ref{NOTcrit} part (iii)], $U_{2}$,
as defined in (\ref{VNwantnormal}), does not converge in distribution
to a
normal. Hence, we can guarantee that, provided $H\geq2/3$, the adjusted
variation in Theorem \ref{adjust} does not converge to a normal. Thus, the
normality of our adjusted estimator in Theorem \ref{thmstat} holds
\textit{if
and only if} $H\in(1/2,2/3)$.
\end{appendix}

%

%
\printaddresses


\begin{thebibliography}{32}

\bibitem{B}
%
\begin{bbook}[vtex]
\bauthor{\bsnm{Beran},~\bfnm{Jan}\binits{J.}}
(\byear{1994}).
\btitle{Statistics for Long-memory Processes}.
\bseries{Monographs on Statistics and Applied Probability}
\bvolume{61}.
\bpublisher{Chapman and Hall}, \baddress{London}.
\bmrnumber{MR1304490}
\end{bbook}
%
\endbibitem

\bibitem{BN}
%
\begin{barticle}[msn]
\bauthor{\bsnm{Breton},~\bfnm{Jean-Christophe}\binits{J.-C.}} \AND
\bauthor{\bsnm{Nourdin},~\bfnm{Ivan}\binits{I.}}
(\byear{2008}).
\btitle{Error bounds on the nonnormal approximation of {H}ermite power
variations of fractional {B}rownian motion}.
\bjournal{Electron. Comm. Probab.}
\bvolume{13}
\bpages{482--493}.
\bmrnumber{MR2447835}
\end{barticle}
%
\endbibitem

\bibitem{BrMa}
%
\begin{barticle}[msn]
\bauthor{\bsnm{Breuer},~\bfnm{Peter}\binits{P.}} \AND
\bauthor{\bsnm{Major},~\bfnm{P{\'e}ter}\binits{P.}}
(\byear{1983}).
\btitle{Central limit theorems for nonlinear functionals of {G}aussian fields}.
\bjournal{J. Multivariate Anal.}
\bvolume{13}
\bpages{425--441}.
\bmrnumber{MR716933}
\end{barticle}
%
\endbibitem

\bibitem{coeur}
%
\begin{barticle}[msn]
\bauthor{\bsnm{Coeurjolly},~\bfnm{Jean-Fran{\c{c}}ois}\binits{J.-F.}}
(\byear{2001}).
\btitle{Estimating the parameters of a fractional {B}rownian motion by discrete
variations of its sample paths}.
\bjournal{Stat. Inference Stoch. Process.}
\bvolume{4}
\bpages{199--227}.
\bmrnumber{MR1856174}
\end{barticle}
%
\endbibitem

\bibitem{DM}
%
\begin{barticle}[msn]
\bauthor{\bsnm{Dobrushin},~\bfnm{R.~L.}\binits{R.~L.}} \AND
\bauthor{\bsnm{Major},~\bfnm{P.}\binits{P.}}
(\byear{1979}).
\btitle{Noncentral limit theorems for nonlinear functionals of {G}aussian
fields}.
\bjournal{Z. Wahrsch. Verw. Gebiete}
\bvolume{50}
\bpages{27--52}.
\bmrnumber{MR550122}
\end{barticle}
%
\endbibitem

\bibitem{Doob}
%
\begin{bbook}[msn]
\bauthor{\bsnm{Doob},~\bfnm{J.~L.}\binits{J.~L.}}
(\byear{1953}).
\btitle{Stochastic Processes}.
\bpublisher{Wiley}, \baddress{New York}.
\bmrnumber{MR0058896}
\end{bbook}
%
\endbibitem

\bibitem{EM}
%
\begin{bbook}[vtex]
\bauthor{\bsnm{Embrechts},~\bfnm{Paul}\binits{P.}} \AND
\bauthor{\bsnm{Maejima},~\bfnm{Makoto}\binits{M.}}
(\byear{2002}).
\btitle{Selfsimilar Processes}.
\bpublisher{Princeton Univ. Press}, \baddress{Princeton, NJ}.
\bmrnumber{MR1920153}
\end{bbook}
%
\endbibitem

\bibitem{GuLe}
%
\begin{barticle}[msn]
\bauthor{\bsnm{Guyon},~\bfnm{Xavier}\binits{X.}} \AND
\bauthor{\bsnm{Le{\'o}n},~\bfnm{Jos{\'e}}\binits{J.}}
(\byear{1989}).
\btitle{Convergence en loi des {$H$}-variations d'un processus gaussien
stationnaire sur {${\bf R}$}}.
\bjournal{Ann. Inst. H. Poincar\'e Probab. Statist.}
\bvolume{25}
\bpages{265--282}.
\bmrnumber{MR1023952}
\end{barticle}
%
\endbibitem

\bibitem{Ha}
%
\begin{barticle}[msn]
\bauthor{\bsnm{Hariz},~\bfnm{Samir Ben}\binits{S. B.}}
(\byear{2002}).
\btitle{Limit theorems for the nonlinear functional of stationary {G}aussian
processes}.
\bjournal{J. Multivariate Anal.}
\bvolume{80}
\bpages{191--216}.
\bmrnumber{MR1889773}
\end{barticle}
%
\endbibitem

\bibitem{HN}
%
\begin{barticle}[msn]
\bauthor{\bsnm{Hu},~\bfnm{Yaozhong}\binits{Y.}} \AND
\bauthor{\bsnm{Nualart},~\bfnm{David}\binits{D.}}
(\byear{2005}).
\btitle{Renormalized self-intersection local time for fractional {B}rownian
motion}.
\bjournal{Ann. Probab.}
\bvolume{33}
\bpages{948--983}.
\bmrnumber{MR2135309}
\end{barticle}
%
\endbibitem

\bibitem{LaIs}
%
\begin{barticle}[msn]
\bauthor{\bsnm{Lang},~\bfnm{Gabriel}\binits{G.}} \AND
\bauthor{\bsnm{Istas},~\bfnm{Jacques}\binits{J.}}
(\byear{1997}).
\btitle{Quadratic variations and estimation of the local {H}\"older
index of a
{G}aussian process}.
\bjournal{Ann. Inst. H. Poincar\'e Probab. Statist.}
\bvolume{33}
\bpages{407--436}.
\bmrnumber{MR1465796}
\end{barticle}
%
\endbibitem

\bibitem{Mand}
%
\begin{barticle}[vtex]
\bauthor{\bsnm{Mandelbrot},~\bfnm{B.}\binits{B.}}
(\byear{1963}).
\btitle{The variation of certain speculative prices}.
\bjournal{J. Bus. Econom. Statist.}
\bvolume{36}
\bpages{392--417}.
\end{barticle}
%
\endbibitem

\bibitem{Leod}
%
\begin{barticle}[vtex]
\bauthor{\bsnm{McLeod},~\bfnm{A.~I.}\binits{A.~I.}} \AND
\bauthor{\bsnm{Kipel},~\bfnm{K.~W.}\binits{K.~W.}}
(\byear{1978}).
\btitle{Preservation of the rescaled adjusted range:
A~reaassement of the Hurst exponent}.
\bjournal{Water Resourc. Res.}
\bvolume{14}
\bpages{491--508}.
\end{barticle}
%
\endbibitem

\bibitem{LeLu}
%
\begin{barticle}[msn]
\bauthor{\bsnm{Le{\'o}n},~\bfnm{Jos{\'e}}\binits{J.}} \AND
\bauthor{\bsnm{Lude{\~n}a},~\bfnm{Carenne}\binits{C.}}
(\byear{2007}).
\btitle{Limits for weighted {$p$}-variations and likewise functionals of
fractional diffusions with drift}.
\bjournal{Stochastic Process. Appl.}
\bvolume{117}
\bpages{271--296}.
\bmrnumber{MR2290877}
\end{barticle}
%
\endbibitem

\bibitem{MaRo}
%
\begin{bmisc}[unstr]
\bauthor{\bsnm{Marcus},~\bfnm{M.~B.}\binits{M.~B.}} \AND
\bauthor{\bsnm{Rosen},~\bfnm{J.}\binits{J.}}
(\byear{2007}).
\bhowpublished{Nonnormal CLTs for functions of the increments of
Gaussian processes with conve increment's variance. Preprint}.
\end{bmisc}
%
\endbibitem

\bibitem{No}
%
\begin{barticle}[unstr]
\bauthor{\bsnm{Nourdin},~\bfnm{I.}\binits{I.}}
(\byear{2008}).
\btitle{Asymptotic behavior of certain weighted quadratic
variation and cubic varitions of fractional Brownian motion}.
\bjournal{Ann. Probab.}
\bvolume{36}
\bpages{2159--2175}.
\end{barticle}
%
\endbibitem

\bibitem{NoNu}
%
\begin{bmisc}[unstr]
\bauthor{\bsnm{Nourdin},~\bfnm{I.}\binits{I.}} \AND
\bauthor{\bsnm{Nualart},~\bfnm{D.}\binits{D.}}
(\byear{2007}).
\bhowpublished{Central limit theorems for multiple Skorohod integrals.
Preprint}.
\end{bmisc}
%
\endbibitem

\bibitem{NoPe}
%
\begin{barticle}[vtex]
\bauthor{\bsnm{Nourdin},~\bfnm{Ivan}\binits{I.}} \AND
\bauthor{\bsnm{Peccati},~\bfnm{Giovanni}\binits{G.}}
(\byear{2008}).
\btitle{Weighted power variations of iterated {B}rownian motion}.
\bjournal{Electron. J. Probab.}
\bvolume{13}
\bpages{1229--1256}.
\bmrnumber{MR2430706}
\end{barticle}
%
\endbibitem

\bibitem{NoPe2}
%
\begin{barticle}[unstr]
\bauthor{\bsnm{Nourdin},~\bfnm{I.}\binits{I.}} \AND
\bauthor{\bsnm{Peccati},~\bfnm{G.}\binits{G.}}
(\byear{2009}).
\btitle{Stein's method on Wiener chaos}.
\bjournal{Probab. Theory Related Fields.}
\bvolume{145}
\bpages{75--118}.
\end{barticle}
%
\endbibitem

\bibitem{NPR}
%
\begin{barticle}[vtex]
\bauthor{\bsnm{Nourdin},~\bfnm{Ivan}\binits{I.}} \AND
\bauthor{\bsnm{R\'{e}veillac},~\bfnm{A.}\binits{G.}}
(\byear{2009}).
\btitle{Multivariate normal approximation using Stein's method and Malliavin
calculus}.
\bjournal{Ann. Inst. H. Poincar\'{e} Probab. Statist.}
To appear.
\end{barticle}
%
\endbibitem

\bibitem{Nbook}
%
\begin{bbook}[vtex]
\bauthor{\bsnm{Nualart},~\bfnm{David}\binits{D.}}
(\byear{2006}).
\btitle{The {M}alliavin Calculus and Related Topics},
\bedition{2nd} ed.
\bpublisher{Springer}, \baddress{Berlin}.
\bmrnumber{MR2200233}
\end{bbook}
%
\endbibitem

\bibitem{NOT}
%
\begin{barticle}[msn]
\bauthor{\bsnm{Nualart},~\bfnm{D.}\binits{D.}} \AND
\bauthor{\bsnm{Ortiz-Latorre},~\bfnm{S.}\binits{S.}}
(\byear{2008}).
\btitle{Central limit theorems for multiple stochastic integrals and
{M}alliavin calculus}.
\bjournal{Stochastic Process. Appl.}
\bvolume{118}
\bpages{614--628}.
\bmrnumber{MR2394845}
\end{barticle}
%
\endbibitem

\bibitem{NP}
%
\begin{barticle}[msn]
\bauthor{\bsnm{Nualart},~\bfnm{David}\binits{D.}} \AND
\bauthor{\bsnm{Peccati},~\bfnm{Giovanni}\binits{G.}}
(\byear{2005}).
\btitle{Central limit theorems for sequences of multiple stochastic integrals}.
\bjournal{Ann. Probab.}
\bvolume{33}
\bpages{177--193}.
\bmrnumber{MR2118863}
\end{barticle}
%
\endbibitem

\bibitem{PT}
%
\begin{bincollection}[vtex]
\bauthor{\bsnm{Peccati},~\bfnm{Giovanni}\binits{G.}} \AND
\bauthor{\bsnm{Tudor},~\bfnm{Ciprian~A.}\binits{C.~A.}}
(\byear{2004}).
\btitle{Gaussian limits for vector-valued multiple stochastic integrals}.
In \bbooktitle{S\'eminaire de {P}robabilit\'es {XXXVIII}}.
\bseries{Lecture Notes in Math.}
\bvolume{1857}
\bpages{247--262}.
\bpublisher{Springer}, \baddress{Berlin}.
\bmrnumber{MR2126978}
\end{bincollection}
%
\endbibitem

\bibitem{SaTa}
%
\begin{bbook}[vtex]
\bauthor{\bsnm{Samorodnitsky},~\bfnm{Gennady}\binits{G.}} \AND
\bauthor{\bsnm{Taqqu},~\bfnm{Murad~S.}\binits{M.~S.}}
(\byear{1994}).
\btitle{Stable Non-{G}aussian Random Variables}.
\bpublisher{Chapman and Hall}, \baddress{London}.
\bmrnumber{MR1280932}
\end{bbook}
%
\endbibitem

\bibitem{Swa}
%
\begin{barticle}[msn]
\bauthor{\bsnm{Swanson},~\bfnm{Jason}\binits{J.}}
(\byear{2007}).
\btitle{Variations of the solution to a stochastic heat equation}.
\bjournal{Ann. Probab.}
\bvolume{35}
\bpages{2122--2159}.
\bmrnumber{MR2353385}
\end{barticle}
%
\endbibitem

\bibitem{Ta1}
%
\begin{barticle}[vtex]
\bauthor{\bsnm{Taqqu},~\bfnm{Murad~S.}\binits{M.~S.}}
(\byear{1975}).
\btitle{Weak convergence to fractional {B}rownian motion and to the
{R}osenblatt process}.
\bjournal{Z. Wahrsch. Verw. Gebiete}
\bvolume{31}
\bpages{287--302}.
\bmrnumber{MR0400329}
\end{barticle}
%
\endbibitem

\bibitem{Ta2}
%
\begin{barticle}[msn]
\bauthor{\bsnm{Taqqu},~\bfnm{Murad~S.}\binits{M.~S.}}
(\byear{1979}).
\btitle{Convergence of integrated processes of arbitrary {H}ermite rank}.
\bjournal{Z.~Wahrsch. Verw. Gebiete}
\bvolume{50}
\bpages{53--83}.
\bmrnumber{MR550123}
\end{barticle}
%
\endbibitem

\bibitem{Tud}
%
\begin{barticle}[msn]
\bauthor{\bsnm{Tudor},~\bfnm{Ciprian~A.}\binits{C.~A.}}
(\byear{2008}).
\btitle{Analysis of the {R}osenblatt process}.
\bjournal{ESAIM Probab. Stat.}
\bvolume{12}
\bpages{230--257}.
\bmrnumber{MR2374640}
\end{barticle}
%
\endbibitem

\bibitem{Ustu}
%
\begin{bbook}[vtex]
\bauthor{\bsnm{{\"U}st{\"u}nel},~\bfnm{Ali~S{\"u}leyman}\binits{A.~S.}}
(\byear{1995}).
\btitle{An Introduction to Analysis on {W}iener Space}.
\bseries{Lecture Notes in Math.}
\bvolume{1610}.
\bpublisher{Springer}, \baddress{Berlin}.
\bmrnumber{MR1439752}
\end{bbook}
%
\endbibitem

\bibitem{WiTaTe}
%
\begin{barticle}[vtex]
\bauthor{\bsnm{Willinger},~\bfnm{W.}\binits{W.}},
\bauthor{\bsnm{Taqqu},~\bfnm{M.}\binits{M.}} \AND
\bauthor{\bsnm{Teverovsky},~\bfnm{V.}\binits{V.}}
(\byear{1999}).
\btitle{Long range dependence and stock returns}.
\bjournal{Finance Stoch.}
\bvolume{3}
\bpages{1--13}.
\end{barticle}
%
\endbibitem

\bibitem{WiTaLeWi}
%
\begin{barticle}[vtex]
\bauthor{\bsnm{Willinger},~\bfnm{W.}\binits{W.}},
\bauthor{\bsnm{Taqqu},~\bfnm{M.}\binits{M.}},
\bauthor{\bsnm{Leland},~\bfnm{W.~E.}\binits{W.~E.}} \AND
\bauthor{\bsnm{Wilson},~\bfnm{D.~V.}\binits{D.~V.}}
(\byear{1995}).
\btitle{Self-similarity in high speed packet traffic: Analysis
and modelisation of ethernet traffic measurements}.
\bjournal{Statist. Sci.}
\bvolume{10}
\bpages{67--85}.
\end{barticle}
%
\endbibitem

\end{thebibliography}
\end{document}